\documentclass[10pt]{article}

\usepackage{amsmath} 
\usepackage{amsthm}
\usepackage{amssymb}
\input epsf

\def\protectbold#1{\protect{\boldmath{$#1$}}}

\theoremstyle{plain}             

\newtheorem{theorem}{Theorem}[section]
\theoremstyle{definition}

\newtheorem{remark}[theorem]{Remark}

\def\eoremark{{\unskip\nobreak\hfil\penalty50	
\hskip2em\hbox{}\nobreak\hfil$\triangle$
\parfillskip=0pt\finalhyphendemerits=0\medbreak}}

\def\Frac#1#2{\frac{{\raise.6ex \hbox{$\displaystyle#1$}}}
{{\lower.6ex\hbox{$\displaystyle#2$}}}}

\def\tfrac#1#2{{{\lower.6ex\hbox{$\scriptstyle#1$}}\over {\raise.7ex\hbox{$\scriptstyle#2$}}}}

\def\dsp{\displaystyle}

\def\eqref#1{(\ref{#1})}

\def\eps{\varepsilon}

\def\wt{\widetilde}

\def\Ai{{{\rm Ai}}}

\def\arccosh{{\rm arccosh}}
\def\arcsinh{{\rm arcsinh}}
\def\arctanh{{\rm arctanh}}

\def\calC{{{\cal C}}}
\def\calL{{{\cal L}}}
\def\bigO{{\cal O}}

\def\ph{{\rm ph}}

\def\binomial#1#2{
\renewcommand{\arraystretch}{1.0}
\left(
\begin{array}{c} 
\hskip-5pt#1\\
\hskip-5pt#2
\end{array}
\hskip-5pt\right)}

\begin{document}
 \title{Asymptotic approximations to the nodes and weights of Gauss--Hermite and Gauss--Laguerre quadratures}

\author{
\hspace*{-0.7cm} A. Gil\\
Departamento de Matem\'atica Aplicada y CC. de la Computaci\'on.\\
ETSI Caminos. Universidad de Cantabria. 39005-Santander, Spain.\\
 \and
J. Segura\\
        Departamento de Matem\'aticas, Estadistica y 
        Computaci\'on,\\
        Universidad de Cantabria, 39005 Santander, Spain.\\
\and
N.M. Temme\\
IAA, 1825 BD 25, Alkmaar, The Netherlands.\footnote{Former address: Centrum Wiskunde \& Informatica (CWI), 
        Science Park 123, 1098 XG Amsterdam,  The Netherlands}\\
}


\maketitle
\begin{abstract}
Asymptotic approximations to the zeros of Hermite and Laguerre 
polynomials are given, together with methods for obtaining the coefficients in the expansions. 
These approximations can be used 
as a standalone method of computation of Gaussian quadratures for high enough degrees, with Gaussian 
weights computed from asymptotic approximations for the orthogonal polynomials. 
We provide numerical evidence showing that for degrees greater than $100$ 
the asymptotic methods are enough for a double precision accuracy computation ($15$-$16$ digits)
of the nodes and weights of the Gauss--Hermite and Gauss--Laguerre quadratures. 
\end{abstract}

\section{Introduction}

As is well known, the nodes $x_i, \, i=1,\ldots, n$ of  Gaussian quadrature rules
are the roots of the (for instance monic) orthogonal polynomial satisfying

\begin{equation}
\label{int:eq1}
\displaystyle\int_a^b x^i p_n(x)w(x)dx = 0,\quad i = 0,\ldots, n-1.
\end{equation}

Among the Gauss quadrature rules, the most popular are those for which the associated
orthogonal polynomials  are the so-called classical orthogonal polynomials, namely: 

\begin{enumerate}
\item Gauss--Hermite:   $w(x) = e^{-x^2}$; $a = -\infty$, $b = +\infty$. Orthogonal polynomials: Hermite polynomials ($H_n(x)$);

\item Gauss--Laguerre: $w(x) = x^{\alpha}e^{-x}$, $\alpha > -1$; $a = 0$, $b = +\infty$. Orthogonal polynomials: Laguerre polynomials ($L_n^{(\alpha)}(x)$);

\item Gauss-Jacobi: $w(x) = (1-x)^{\alpha}(1+x)^{\beta}$, $\alpha,\,\beta >-1$; $a = -1$, $b = 1$. Orthogonal
polynomials: Jacobi polynomials ($P_n^{(\alpha,\beta)}(x)$).
\end{enumerate} 

The weights for
the $n$-point Gauss quadrature based on the nodes $\{x_i\}_{i=1}^{n}$  can be written 
in terms of the derivatives of the orthogonal polynomials at the nodes as follows:

\begin{enumerate}
\item Gauss--Hermite:  

\begin{equation}
\label{wgh}
w_i=\Frac{\sqrt{\pi}2^{n+1}n! }{ [H_n^{\prime}(x_i)]^2},
\end{equation}

\item Gauss--Laguerre:

\begin{equation}
\label{wgl}
w_i = \Frac{\Gamma (n+\alpha+1)}{n! x_i [L_{n} ^{(\alpha)\prime}(x_i)]^2},
\end{equation}

\item Gauss-Jacobi: 

\begin{equation}
\label{wgj}
w_i  =\Frac{M_{n,\alpha,\beta}}{(1-x_i^2) [P_{n} ^{(\alpha ,\beta) \prime}(x_i)]^2},
\end{equation}
where 
$$
M_{n,\alpha,\beta}=2^{\alpha+\beta+1}\Frac{\Gamma (n+\alpha+1)\Gamma (n+\beta+1)}{n!\Gamma (n+\alpha+\beta+1 )}.
$$
\end{enumerate}

Iterative algorithms are interesting methods of computation of Gaussian nodes and weights, very clearly 
outperforming matrix methods (Golub-Welsch \cite{Golub:1969:COG}) for high degrees. They are based 
on the computation of the roots of the orthogonal polynomial by an iterative method and the
subsequent computation of the weights by using function relations like those in 
Eqs. (\ref{wgh})-(\ref{wgj}).
Most iterative methods for the computation of the Gaussian 
nodes (with the exception of \cite{Segura:2010:RCO}) 
require accurate enough first approximations in order
to ensure the convergence of the iterative method (typically the Newton method); for two recent
examples, see \cite{Hale:2013:SJSC,Town:2016:IMA}. An alternative approach 
\cite{Glaser:2007:AFA}, although less efficient for
high degrees than iterative methods with asymptotic first approximations \cite{Hale:2013:SJSC,Town:2016:IMA}, consists in guessing these first approximations
by integrating a Prufer-transformed ODE with a Runge-Kutta method, and then refining these
guesses by the Newton method (however, asymptotic approximations were also used in this reference for
the particular case of Gauss-Legendre quadrature). More recently, 
non-iterative methods based on asymptotic approximations 
for the computation of Gauss-Legendre nodes and weights were developed in \cite{Bogaert:2014:IFC}, which
were shown to outperform iterative approaches. 

In this paper, our aim is to provide asymptotic approximations
 for the accurate computation of the nodes and weights of Gauss--Hermite and Gauss--Laguerre
quadrature. These approximations provide a fast and accurate method of computation which 
can be used for arbitrarily large degree, but which
also provide accurate results for not so large degrees ($n\ge 100$). 
The methods are able to compute both the nodes and the weights with nearly double precision accuracy, improving the accuracy of the available fixed precision iterative methods.

As we will discuss in
a subsequent paper, a fully non-iterative approach is also possible for the 
case of Gauss-Jacobi quadrature \cite{Gil:2017:NCO}, similarly as was shown for the
particular case of Legendre polynomials \cite{Bogaert:2014:IFC}.

\section{Hermite polynomials}\label{sec:Hermpol}

In \cite{Town:2016:IMA}  first estimates of the zeros of Hermite polynomials
are based on work of Tricomi for the middle zeros; these first guesses follow from expansions in terms of elementary functions. For the remaining zeros near the positive endpoint $\sqrt{2n+1}$ of the zeros interval the first estimates are taken from the work of Gatteschi, and are in terms of the zeros of the Airy functions.

In this section we give an expansion of the zeros based on the asymptotic expansion  in terms of Airy functions. The expansion can be used for all positive zeros, however, the approximations are less accurate for the small zeros. For these we give an approximation
 based on an asymptotic expansion in terms of elementary functions. We start discussing this expansion.

 \subsection{Expansions in terms of elementary functions}\label{sec:element}
An expansion in terms of elementary functions for the Hermite polynomials is given in \cite[\S18.15(v)]{koorn:2010:OPS} with a limited number of coefficients. However, we prefer an expansion for the parabolic cylinder function derived in \cite{Olver:1959:UAE}; these results are summarized in \cite[\S12.10(iv)]{Temme:2010:PCF} and \cite[\S30.2.3]{Temme:2015:AMI}.

The relation between the parabolic cylinder function $U(a,z)$ and the Hermite polynomial $H_n(z)$ is
\begin{equation}\label{eq:herelem01}
U\left(-n-\tfrac12,z\right)=2^{-n/2}e^{-\frac14z^2}H_n\left(z/\sqrt2\right),\quad n=0,1,2,\ldots.
\end{equation}
We use the notations
\begin{equation}\label{eq:herelem02}
\mu=\sqrt{2n+1},\quad t=x/\mu,\quad \eta(t)=\tfrac12\arccos t-\tfrac12 t\sqrt{1-t^2},
\end{equation}
and we have the asymptotic representation
\begin{equation}\label{eq:herelem03}
\begin{array}{@{}r@{\;}c@{\;}l@{}}
H_n(x)&=&
\dsp{\frac{2^{\frac12n+1}e^{\frac12x^2}g(\mu)}{(1-t^2)^{\frac14}}}\ \times\\[8pt]
&&\left(\cos\left(\mu^2\eta-\tfrac14\pi\right)
{{\cal A}}_\mu(t)
-\sin\left(\mu^2\eta-\tfrac14\pi\right)
{{\cal B}}_\mu(t)\right),
\end{array}
\end{equation}
with expansions
\begin{equation}\label{eq:herelem04}
\begin{array}{ll}
\dsp{
{{\cal A}}_\mu(t)\sim\sum_{s=0}^\infty\frac{(-1)^su_{2s}(t)}{(1-t^2)^{3s}\mu^{4s}},\quad
{{\cal B}}_\mu(t)\sim\sum_{s=0}^\infty\frac{(-1)^su_{2s+1}(t)}{(1-t^2)^{3s+\frac32}\mu^{4s+2}},}
\end{array}
\end{equation}
uniformly for $-1+\delta\le t\le 1-\delta$, where $\delta$ is an arbitrary small positive number.

The first few coefficients are
\begin{equation}\label{eq:herelem05}
u_0(t)=1,\quad u_1(t)=\frac{t(t^2-6)}{24}, \quad
u_2(t)=\frac{-9t^4+249t^2+145}{1152},
\end{equation}
and more $u_s(t)$ follow from  the recurrence relations
\begin{equation}\label{eq:herelem06}
\begin{array}{ll}
(t^2-1)u'_s(t)-3stu_s(t)=r_{s-1}(t),\\[8pt]
8r_s(t)=(3t^2+2)u_s(t)-12(s+1)tr_{s-1}(t)+4(t^2-1)r'_{s-1}(t).
\end{array}
\end{equation}

The quantity $g(\mu)$ is only known in the form of an asymptotic expansion
\begin{equation}\label{eq:herelem07}
g(\mu)\sim h(\mu)\left(1+\tfrac12\sum_{k=0}^\infty \frac{\gamma_k}{\left(\frac12\mu^2\right)^k}\right),
\end{equation}
where the coefficients $\gamma_k$ are defined by
\begin{equation}\label{eq:herelem08}
\Gamma\left(\tfrac12+z\right)\sim \sqrt{2\pi}\, e^{-z}\,z^z\,\sum_{k=0}^\infty \frac{\gamma_k}{z^k}, \quad z\to\infty.
\end{equation}
The first ones are
\begin{equation}\label{eq:herelem09}
\gamma_0=1,\quad \gamma_1 =- \tfrac{1}{24}, \quad \gamma_2=\tfrac{1}{1152},\quad 
\gamma_3 = \tfrac{1003}{414720}, \quad \gamma_4=-\tfrac{4027}{39813120}.
\end{equation}
For $h(\mu)$ we have
\begin{equation}\label{eq:herelem10}
h(\mu)=2^{-\frac14\mu^2-\frac14}e^{-\frac14\mu^2}\mu^{\frac12\mu^2-\tfrac12}=2^{-\frac12}\left(n+\tfrac12\right)^{\frac12n}e^{-\frac12n-\frac14}.
\end{equation}

\subsubsection{Expansions of the zeros}\label{sec:zerelem}
Next we discuss expansions for the zeros of $H_n(x)$, $x_k$, $1\le k \le n$ ($x_1<x_2<\cdots<x_n)$. We introduce a function $W(\eta)$ (see \eqref{eq:herelem03})
\begin{equation}\label{eq:herelem11}
W(\eta)=\cos\left(\mu^2\eta-\tfrac14\pi\right){{\cal A}}_\mu(t)
-\sin\left(\mu^2\eta-\tfrac14\pi\right)
{{\cal B}}_\mu(t),
\end{equation}
and try to solve the equation $W(\eta)=0$ for large values of $n$. We define a first approximation $\eta_0$ such that the cosine term vanishes and $\eta_0$ and the corresponding $t$ and $x$-values are  
(in first-order approximation) related to 
a zero of $H_n(x)$. 

The small zeros are around $x=0$ and $t=0$, that is, for $\eta$ near $\eta(0)=\frac14\pi$.  We define
\begin{equation}\label{eq:herelem12}
\eta_0=\frac{n-k+\frac{3}{4}}{\mu^2}\pi,\quad k=1,2,\ldots,n.
\end{equation}
In this way, $\cos\left(\mu^2\eta_0-\frac14\pi\right)=0$, and this choice of $\eta_0$ follows from the location of the zeros of the cosine function and those of $H_n(x)$. 
Observe that, when $n$ is odd and $k=\frac12(n+1)$, that is, $x_k=0$, it follows that $\eta_0=\frac14\pi$. If $\eta=\frac14\pi$ we have $t=0$ and $x=0$.

We assume that the equation $W(\eta)=0$ has a solution $\eta$ that can be expanded in the form
\begin{equation}\label{eq:herelem13}
\eta=\eta_0+\eps,\quad \eps\sim \frac{\eta_1}{\mu^2}+\frac{\eta_2}{\mu^4}+\frac{\eta_3}{\mu^6}+\frac{\eta_4}{\mu^8}+\ldots,
\end{equation}
and consider the Taylor expansion and equation
\begin{equation}\label{eq:herelem14}
W(\eta)+\frac{\eps}{1!}\frac{d}{d\eta}W(\eta)+\frac{\eps^2}{2!}\frac{d^2}{d\eta^2}W(\eta)+\frac{\eps^3}{3!}\frac{d^2}{d\eta^3}W(\eta)+\ldots=0,
\end{equation}
where $W(\eta)$ and its derivatives are taken at $\eta=\eta_0$. Because the expansions in \eqref{eq:herelem04} are in terms of $t$,  we need $dt/d\eta=-1/\sqrt{1-t^2}$.

When we have found $\eta$, the corresponding $t$-value is obtained by inverting the relation for $\eta(t)$ in \eqref{eq:herelem02}. For this purpose we use the expansion
\begin{equation}\label{eq:herelem15}
\begin{array}{l}
t=-\wt\eta-\frac16\wt\eta^3-\frac{13}{120}\wt\eta^5-\frac{493}{5040}\wt\eta^7+\cdots,\\
\wt\eta=\eta-\tfrac14\pi=-\frac{1}{2}\arcsin t-\frac{1}{2}t\sqrt{1-t^2}
\end{array}
\end{equation}

It is also possible to invert the relation (\ref{eq:herelem15}) by using an iterative method. 
For this purpose it is convenient to write $t=\sin\frac12\theta$. Then the equation to be solved for $\theta\in(-\pi,\pi)$ reads 
\begin{equation}\label{eq:herelem16}
4\wt\eta+\theta+\sin\theta=0.
\end{equation}
A Newton or related procedure can be used to solve this equation, but in our algorithms we prefer
to use the series shown in (\ref{eq:herelem15}), which is faster (and of more restricted applicability, but
sufficient for our purposes).

After a few symbolic manipulations we find that $\eta_{2k+1}=0$, $k=0,1,2,\ldots$, and that the first nonzero coefficients are
\begin{equation}\label{eq:herelem17}
\begin{array}{@{}r@{\;}c@{\;}l@{}}
\eta_2&=&\dsp{ -\frac{t\left(t^2-6\right)}{24\left(1-t^2\right)^{\frac32}},}\\[8pt]
\eta_4&=&\dsp{ -\frac{t\left(56t^8-252t^6 +351t^4+ 2340t^2+3780\right)}{5760\left(1-t^2\right)^{\frac92}},}\\[8pt]
\eta_6&=&-t\bigl(3968t^{14}-29760t^{12}+95544t^{10}-173232t^8+231237t^6\ - \\[8pt]
&&1890882t^4-6068580t^2-1690920\bigr)/\left(322560(1-t^2\right)^{\frac{15}{2}}).
\end{array}
\end{equation}
Because we have a recurrence relation for the coefficients $u_s(t)$ in \eqref{eq:herelem06}, it is quite easy to generate many $u_s(t)$ and also much more coefficients $\eta_j$ than given in \eqref{eq:herelem17}.

\paragraph{Algorithm}
For the computation of the approximations of the zeros $x_k$ we summarize the procedure as follows.
\begin{enumerate}
\item
To approximate the zero $x_{k}$, compute the starting value $\eta_0 $, given in \eqref{eq:herelem12}.

\item
Compute the corresponding $t$-value from \eqref{eq:herelem15} (with $\eta=\eta_0$). 

\item
With these values $\eta_0$ and $t$, compute the coefficients $\eta_k$ in \eqref{eq:herelem17}. 

\item
Next, compute $\eta$  from \eqref{eq:herelem13}. 

\item
Then the better value of $t$ again follows from \eqref{eq:herelem15}. 

\item
Finally,  the approximation for the requested zero is $x_k\sim \mu t$, see \eqref{eq:herelem02}. 

\end{enumerate}

 \subsection{Expansions in terms of Airy functions}\label{sec:airy}
For the large zeros we shall use the Airy-type expansion of  the Hermite polynomials. We write (see \cite[Section 12.10(vii)]{Temme:2010:PCF}) 
\begin{equation}\label{eq:HermAiry01}
H_n(x)=\sqrt{\pi}\,2^{\frac12n+1}\mu^{\frac13}\chi(\zeta)e^{\frac12x^2}g(\mu)\Bigl(\Ai\left(\mu^{\frac43}\zeta\right)A(\zeta)\ +\mu^{-\frac83}\Ai^{\prime}\left(\mu^{\frac43}\zeta\right)B(\zeta)\Bigr),
\end{equation}
with expansions
\begin{equation}\label{eq:HermAiry02}
A(\zeta)\sim\sum_{s=0}^\infty\frac{A_s(\zeta)}{\mu^{4s}}, \quad B(\zeta)\sim \sum_{j=0}^\infty \frac{B_s(\zeta)}{\mu^{4s}},\quad \mu\to\infty,
\end{equation}
where $\mu=\sqrt{2n+1}$,  $t=x/\mu$, and $g(\mu)$ is the function with asymptotic expansion given in \eqref{eq:herelem07}. 
For $\zeta$  we have the definition
\begin{equation}\label{eq:HermAiry03}
\begin{array}{l}
\frac23\zeta^{\frac32}=\frac12t\sqrt{t^2-1}-\frac12\arccosh\,t,\,t\ge1,\\
\\
\frac23(-\zeta)^{\frac32}=\eta(t),\,-1<t\le1,
\end{array}
\end{equation}  
where $\eta(t)$ is defined in \eqref{eq:herelem02}; $\chi(\zeta)$ is defined by
\begin{equation}\label{eq:HermAiry04}
\chi(\zeta)=\left(\frac{\zeta}{t^2-1}\right)^{\frac14}.
\end{equation}  

The variable $\zeta$ is analytic in a neighborhood of $t=1$. We have the differential equation
\begin{equation}\label{eq:HermAiry05}
\zeta\left(\frac{d\zeta}{dt}\right)^2=t^2-1,
\end{equation}
and we have the following expansions in powers of $t-1$ and $\zeta$ the expansions
\begin{equation}\label{eq:HermAiry06}
\begin{array}{@{}r@{\;}c@{\;}l@{}l@{\,}}
2^{-\frac13}\zeta&=&\dsp{(t-1)+\tfrac1{10}(t-1)^2-\tfrac2{175}(t-1)^3+\cdots},\\[8pt]
t&=&\dsp{1+\wt\zeta-\tfrac{1}{10}\wt\zeta^2+\tfrac{11}{350}\wt\zeta^3+\cdots,\quad  \wt\zeta=2^{-\frac13}\zeta}.
\end{array}
\end{equation}
The relation between $t$ and $\zeta$ is singular at $t=-1$, $\zeta(-1)=-(3\pi/4)^{2/3}=-1.770\cdots$, and the series in the second line converges for $\vert\zeta\vert<1.770\cdots$. 

The coefficients are given by
\begin{equation}\label{eq:HermAiry07}
\begin{array}{@{}r@{\;}c@{\;}l@{}l@{\,}}
A_s(\zeta)&=& \dsp{\zeta^{-3s}\sum_{m=0}^{2s}\beta_m\,\left(\chi(\zeta)\right)^{6(2s-m)}u_{2s-m}(t)}, \\[8pt]
B_s(\zeta)&=& \dsp{-\zeta^{-3s-2}\sum_{m=0}^{2s+1}\alpha_m\,\left(\chi(\zeta)\right)^{6(2s-m+1)}u_{2s-m+1}(t)}, 
\end{array}
\end{equation}
where $u_s(t)$ are as in \S\ref{sec:element}, and
\begin{equation}\label{eq:HermAiry08}
\begin{array}{@{}r@{\;}c@{\;}l@{}l@{\,}}
\alpha_m&=& \dsp{\frac{(2m+1)(2m+3)\cdots(6m-1)}{m!\,(144)^m}}, \quad \alpha_0=1,\\[8pt]
\beta_m&=&\dsp{-\frac{6m+1}{6m-1}\alpha_m. }
\end{array}
\end{equation}
A recursion for $\alpha_m$ reads
\begin{equation}\label{eq:HermAiry09}
\alpha_{m+1}=\alpha_m\frac{(6m+5)(6m+3)(6m+1)}{144(m+1)(2m+1)}, \quad m=0,1,2,\ldots\,.
\end{equation}

The first few coefficients of the expansions in \eqref{eq:HermAiry02} are given by:
\begin{equation}\label{eq:HermAiry10}
\begin{array}{@{}r@{\;}c@{\;}l@{}l@{\,}}
A_0(\zeta)&=& 1,\quad \dsp{B_0(\zeta)=-\frac{48\chi^6u_1(t)+5}{48\,\zeta^2},}  \\[8pt]
A_1(\zeta)&=& \dsp{\frac{4608\chi^{12}u_2(t)-672\chi^6u_1(t)-455}{4608\,\zeta^3},}  \\[8pt]
B_1(\zeta)&=& \dsp{-\frac{663552\chi^{18}u_3(t)+69120\chi^{12}u_2(t)+55440\chi^6u_1(t)+85085}{6635528\,\zeta^5}.}  
\end{array}
\end{equation}
Here $\chi=\chi(\zeta)$ is given by \eqref{eq:HermAiry04}. To avoid numerical cancellations when $\zeta$ is small in the above representations, 
we can expand the coefficients, which are analytic at $\zeta=0$, in powers of $\zeta$.

\subsubsection{Expansions of the zeros}\label{sec:Hermsimzer}

An expansion for the zeros is obtained as follows. First we determine the zeros in terms of $\zeta$. 

For the first-order approximation of a zero $x_{n-k+1}$ of $H_n(x)$ we  compute $\zeta_0=\mu^{-\frac43}a_k$, where $a_k$ is a zero of the Airy function $\Ai(x)$. Because of the symmetry of the Hermite polynomial, we assume that $1\le k \le \lfloor \frac12n\rfloor$. 

We introduce an expansion of  $\zeta$ corresponding to the zero of $H_n(x)$ by writing
\begin{equation}\label{eq:HermAiry11}
\zeta=\zeta_0+\eps,\quad \eps\sim \frac{\zeta_1}{\mu^4}+\frac{\zeta_2}{\mu^8}+\ldots,
\end{equation}
and we try to obtain the $\zeta_j, j\ge1$. We introduce a function $W(\zeta)$ by writing (see \eqref{eq:HermAiry01})
\begin{equation}\label{eq:HermAiry12}
 W(\zeta)=\Ai\left(\mu^{\frac43}\zeta\right)A(\zeta)\ +\mu^{-\frac83}\Ai^{\prime}\left(\mu^{\frac43}\zeta\right)B(\zeta),
\end{equation}
and expand $W(\zeta)$ at $\zeta=\zeta_0$, writing $\zeta=\zeta_0+\eps$, which gives 
\begin{equation}\label{eq:HermAiry13}
W(\zeta_0)+\frac{\eps}{1!}W^\prime(\zeta_0)+ \frac{\eps^2}{2!}W^{\prime\prime}(\zeta_0)+\ldots = 0.
\end{equation}
In this equation we substitute the expansion given in \eqref{eq:HermAiry11} and those in \eqref{eq:HermAiry02}, compare equal powers of $\mu$ and obtain the first few coefficients
\begin{equation}\label{eq:HermAiry14}
\begin{array}{@{}r@{\;}c@{\;}l@{}l@{\,}}
\zeta_ 1&=& -B_0(\zeta_0), \\[8pt]
\zeta_2&=& -\frac{1}{3}\left(3B_1(\zeta_0)-3B_0(\zeta_0) A_1(\zeta_0)-3 B_0(\zeta_0) B_0^\prime(\zeta_0)+\zeta_0 B_0(\zeta_0)^3\right),
\end{array}
\end{equation}
where the derivative is with respect to $\zeta$ and the coefficients are given in \eqref{eq:HermAiry07}. 

To obtain the derivative of $B_0(\zeta)$ we need
\begin{equation}\label{eq:HermAiry15}
\frac{dt}{d\zeta}=\chi^2(\zeta),\quad \frac{d\chi}{d\zeta}=\frac{1-2t\chi^6(\zeta)}{4\zeta}\chi(\zeta),
\end{equation}
which  follow from  \eqref{eq:HermAiry04} and \eqref{eq:HermAiry05}. This gives
\begin{equation}\label{eq:HermAiry16}
\frac{d}{d\zeta}B_0(\zeta)=\frac{\chi^6t^3+6\chi^{12}t^4-6t\chi^6-36t^2\chi^{12}-6\chi^8\zeta t^2+12\chi^8\zeta+10}{48\zeta^3}.
\end{equation}
For small values of $\zeta$ we have expansions of the form
\begin{equation}\label{eq:HermAiry17}
\begin{array}{@{}r@{\;}c@{\;}l@{}l@{\,}}
\zeta_ 1&=& 2^{\frac13}\left(\frac{9}{280}-\frac{7}{450}\wt\zeta+\frac{1359}{134750}\wt\zeta^2+\ldots\right), \quad \wt\zeta=2^{-\frac13}\zeta,\\[8pt]
\zeta_2&=&2^{\frac13}\left(-\frac{1539}{130000}+\frac{1550191}{138915000}\wt\zeta-\frac{ 193351}{16362500}\wt\zeta^2+\ldots\right).
\end{array}
\end{equation}

\paragraph{Algorithm}
When we have obtained a value $\zeta$ that corresponds to a zero of the Hermite polynomial, the corresponding $t$-value should be obtained from the second equation in \eqref{eq:HermAiry03}. This equation has to be solved by a numerical procedure. A first estimate, when $\zeta$ is small, can be obtained from the second line in \eqref{eq:HermAiry06}, and more terms of that expansion can easily be obtained by a symbolic package.

For an iterative procedure  it is convenient to substitute  $t=\cos\frac12\theta$, with $\theta\in[0,2\pi)$. Then the equation to be solved for $\theta$ reads $\frac{8}{3}(-\zeta)^{\frac{3}{2}}=\theta-\sin\theta$ 
and we can use, for instance, the Newton method for this purpose. However,
in our algorithms we prefer to invert using enough terms in (\ref{eq:HermAiry06}), which is a faster method.

We proceed as follows for computing approximations for the zeros.
\begin{enumerate}

\item
To approximate the zero $x_{n-k+1}$, define the starting value $\zeta_0=\mu^{-\frac43}a_k$, $1\le k \le \frac12 n$, where $a_k$ is a zero of the Airy function.

\item
Compute $t$ from the second line of \eqref{eq:HermAiry06}. 

\item
With these values $\zeta_0$ and $t$, compute the coefficients $\zeta_j$ in \eqref{eq:HermAiry14} and $\chi(\zeta_0)$ from   \eqref{eq:HermAiry04}. 

\item
Next, compute $\zeta$  from \eqref{eq:HermAiry11}. 

\item
Then the better value of $t$ again follows from the second line of \eqref{eq:HermAiry06}. 

\item
Finally, $x_{n-k+1}\sim t\mu$.

\end{enumerate}

\subsection{Numerical performance of the expansions}

The approximation  (\ref{eq:herelem13})  (obtained from the expansion in terms of elementary functions) is accurate
for large $n$ and particularly for the small zeros. As a first numerical example of 
the accuracy, even for quite small $n$, we take $n=11$, $k=7$ (the smallest positive zero). Then, $\eta_0=0.648807$ and  the corresponding $t$ and $x$-values are 
$0.137021$ and $0.657129$. The seventh zero of $H_{11}(x)$ is $0.656810\ldots$, and the relative error is $0.00048$. 
With the shown coefficients in \eqref{eq:herelem17}  we obtain 
$\eta=0.6488732440401913$ and $x =0 .6568095658827670$, with relative error $ 1.52\times10^{-9}$. The computations are done with Maple, with Digits\,=\,16.
With $n=51$ and $k=27$ (the smallest positive zero), the relative error becomes $10^{-15}$.

The expansions in \eqref{eq:herelem04} are uniformly valid for $-1+\delta\le t\le 1-\delta$, where $\delta$ is an arbitrary small positive number. Hence, for the large zeros this method is not reliable, and we need to restrict the number of zeros that we can compute. For example, we can request that $\vert t\vert \le \frac12$, the corresponding $\eta$-value satisfies $\vert\eta-\frac14\pi\vert\le \frac1{12}\pi+\frac18\sqrt{3}=0.478$. When we use the first estimate $\eta_0$ given in \eqref{eq:herelem12} in the equation $\vert\eta_0-\frac14\pi\vert\le 0.478$, we find for $k$ the bound $\vert \frac12n-k\vert\lesssim  \frac{0.478}{\pi}(2n+1)=0.304 n+0.152$. This says that roughly $0.3n$ of the positive zeros can be computed by using the asymptotic approximations of \S\ref{sec:element}, when we request  $\vert t\vert \le \frac12$. In practice, as we will see later, the expansions in terms of elementary functions
can be used for larger values of $|t|$ and when they are accurate, they are preferable
to the expansions in terms of Airy functions because the algorithm is faster.

More extensive tests of the expansions have been performed using finite
precision implementations coded in Fortran 90. In these implementations
only non-iterative methods (power series) are used for the inversion of the variables.

\begin{figure}
\begin{center}
\hspace*{-2cm}
\epsfxsize=16cm \epsfbox{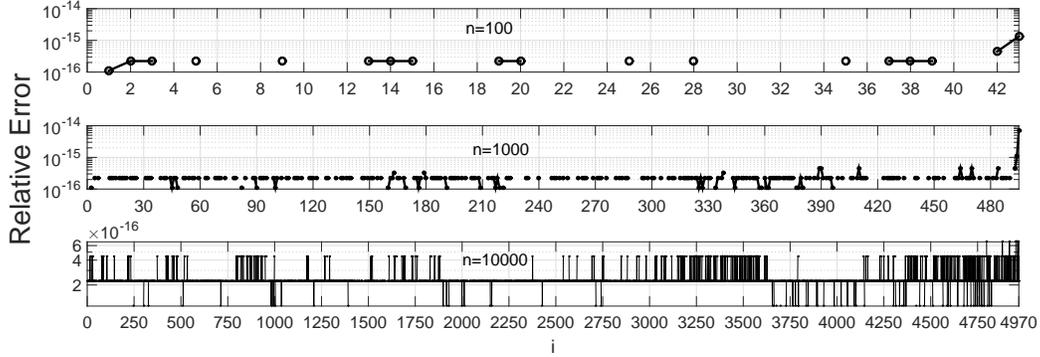}
\caption{ Relative accuracy obtained with the asymptotic expansion (\ref{eq:herelem13}) for computing the positive zeros of $H_n(x)$
for $n=100,\,1000,\,10000$. 
The
label $i$ in the abscissa represents the order of the zero (starting from $i=1$ for the smallest positive zero). 
The points not shown in the plots correspond to values with all digits correct in double precision accuracy. 
\label{Fig1}}
\end{center}
\end{figure}

Figure \ref{Fig1} shows the performance of the expansion in terms of elementary 
functions. In this figure, the relative accuracy obtained for
 computing the positive zeros of $H_n(x)$ for $n=100,\,1000,\,10000$ is plotted.  The
label $i$ in the abscissa represents the order of the zero (starting from $i=1$ for the smallest positive zero). 
The algorithm for testing the accuracy of the zeros has been implemented in finite precision arithmetic using the
first 6 non-zero terms in the expansion. We compare the asymptotic expansions 
against an extended precision accuracy (close to $32$ digits) 
iterative algorithm which uses the global fixed point method of 
\cite{Segura:2010:RCO}, with orthogonal
polynomials computed by local Taylor series.

As can be seen,
a very large number of the zeros  for the three values of $n$ tested
can be computed with the expansion with a relative accuracy near full double precision.
Actually, the points not shown in the plot correspond to values with all digits correct in double precision accuracy. However, the expansion fails for the largest zeros, as expected.

\begin{figure}
\begin{center}
\hspace*{-1.5cm}
\epsfxsize=14cm \epsfbox{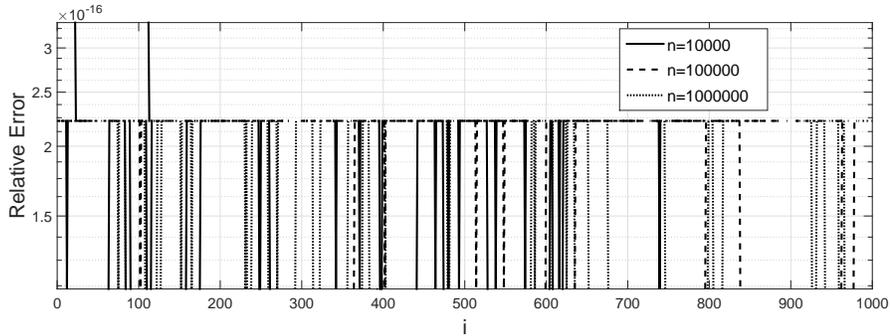}
\caption{ Relative accuracy obtained with the asymptotic expansion (\ref{eq:HermAiry11}) for computing the last $1000$ positive zeros of $H_n(x)$
for $n=10000,\,100000,\,1000000$.
\label{Fig2}}
\end{center}
\end{figure}

As for the asymptotic expansion in terms of the
zeros of Airy functions (\ref{eq:HermAiry11}), the situation is the reverse: the further we are from the turning point at $t=1$ ($\zeta=0$), the larger 
the relative errors become. Therefore, for $n$ fixed the maximum errors in the computation are obtained for the small zeros.
 For example, using Maple with Digits\,=\,16, we take $n=11$ and $6$ coefficients in \eqref{eq:HermAiry11}. Then we have for the zero $x_6$  at the origin $\zeta_0=\mu^{-\frac43}a_6=-1.115618210110694$, $t=  -0.1668495251592333\times 10^{-3}$, and the better values 
$\zeta=-1.115460237225190$ and $t=1.746192313216916 \times 10^{-13}$. This gives 
$x_6\doteq 8.374444141492045\times 10^{-13}$ 
and for the largest zero $x_{11}$ the relative accuracy is $10^{-15}$. A test of the expansion for very large values of $n$ using a finite precision arithmetic implementation is shown in 
Figure \ref{Fig2}. In this figure, we show the relative accuracy  obtained with the asymptotic expansion (\ref{eq:HermAiry11}) for computing the largest
$1000$ positive zeros of $H_n(x)$ for $n=10000,\,100000,\,1000000$. As can be seen, an accuracy near $10^{-16}$ can be obtained
in all cases. The zeros $a_k$ of the Airy function have been computed using 
$a_{k}=-T\left(\tfrac{3}{8}\pi(4k-1)\right)$, where $T(t)$  has the Poincar\'e's expansion (see \cite[\S9.9(iv)]{Olver:2010:ARF})

\begin{equation}\label{AiryPoincare}
T(t)\sim t^{2/3}\left(1+\frac{5}{48}t^{-2}-\frac{5}{36}t^{-4}+\frac{77125}{82944}t^{-6}-\frac{108056875}{6967296}t^{-8}+\cdots\right).
\end{equation}  

This expansion is valid for moderate/large values of $k$. In our implementation we use pre-computed
values for the first 10 zeros of the Airy function and the Poincar\'e's expansion for the rest.  

\begin{figure}
\begin{center}
\hspace*{-1.5cm}
\epsfxsize=14cm \epsfbox{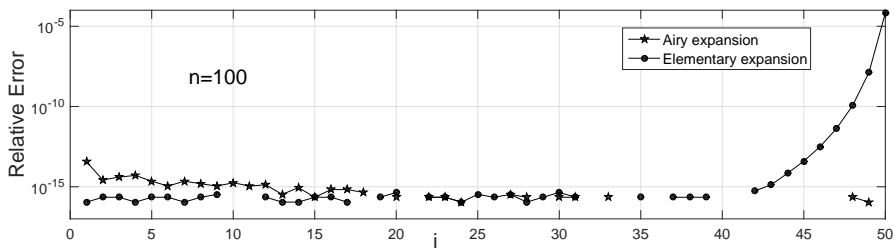}
\caption{ Relative accuracy obtained with the asymptotic expansions (\ref{eq:herelem13})  and (\ref{eq:HermAiry11}) 
for computing the positive zeros of $H_n(x)$ for $n=100$.
The points not shown in the plot correspond to values with all digits correct in double precision accuracy. 
\label{Fig2b}}
\end{center}
\end{figure}

The accuracy of the two expansions  (\ref{eq:herelem13})  and (\ref{eq:HermAiry11})  for approximating 
the zeros of Hermite polynomials for $n=100$ is compared
in Figure \ref{Fig2b}. As can be seen, the combined use of both expansions allow the computation of all the
zeros with a double precision accuracy of $15$-$16$ digits.

In Table \ref{tab:tab02} we illustrate the efficiency of the expansions for approximating the zeros of Hermite polynomials  
for $n=100,\,10000$. In particular, the first~$0.6n$ zeros of the Hermite polynomials have been computed
with the asymptotic expansion in terms of elementary functions and the last $0.4n$ zeros with
the asymptotic expansion in terms of the zeros of Airy functions. With this splitting and by taking enough terms, it is possible to use the series  (\ref{eq:herelem15}) and (\ref{eq:HermAiry06}) for computing the $t$-values
in the expansions instead of using an iterative method for solving the non-linear equations.  
In the table we show average
CPU times (obtained using an Intel Core i54310U 2.6GHz processor under Windows) per 
node. The second column shows the CPU times when the number of terms required (no more than five or six
depending on the expansion)
for a double precision accuracy for the zeros is considered,
while the first column shows the CPU times for only two terms.  For $n=10000$ this is the number
of terms needed in the expansions to obtain double precision accuracy.  
For $n=100$ we observe
that there is not much difference in speed between the more simple ($2$ terms) and the 
more accurate approximation; this favors the use of accurate asymptotic
approximations with no ulterior iterative refinements.
 The table also shows that the computation of the expansion 
in terms of elementary functions is more efficient  than the expansion in terms of zeros of Airy functions
although for $n=10000$ the difference in speed is not very significant.

\renewcommand{\arraystretch}{1.2}
\begin{table}
\caption{Hermite expansions: average CPU times per node. 
The algorithms have been implemented in Fortran 90.
The nodes
are computed with $15-16$ digits accuracy.
\label{tab:tab02}}
$$
\begin{array}{rcc}
{\rm Expansion} & 2\  {\rm terms} \ & \ge 2\ {\rm terms} \quad\\
\hline
  & n=100 & \\
\hline
{\rm Elementary}  & 0.56\,\mu s    &     0.65\,\mu s\\
{\rm Airy}  &   1.0\,\mu s  &   1.24\,\mu s  \\
\hline
& n=10000 & \\
\hline
{\rm Elementary}  &   0.5\,\mu s  &   0.5\,\mu s   \\
{\rm Airy}  &   0.78\,\mu s  &   0.78\,\mu s  \\
\end{array}
$$
\end{table}
\renewcommand{\arraystretch}{1.0}

Once the nodes (the zeros of $H_n(x)$) of the Gauss--Hermite quadrature have been computed, approximations to the weights given in (\ref{wgh}) 
can be also obtained by using the asymptotic results in \S\ref{sec:element} (elementary functions) 
and \S\ref{sec:airy} (Airy functions).

For the computation of the weights, one needs to be careful in order to avoid 
overflows in the computation both as a function of $n$ and as a function of the
values of the nodes. With respect to the dependence on $n$, we observe that the
large factor $2^n n!$ in (\ref{wgh}) can be cancelled out by the factors in
front of the expansions (\ref{eq:herelem03}) and (\ref{eq:HermAiry01}). This is as 
expected because using the first approximations from the elementary asymptotic
expansions as $n\rightarrow \infty$ we obtain the estimate for the weights:
\begin{equation}
w_i\sim \Frac{\pi}{\sqrt{2n}}e^{-x_i^2}.
\end{equation}
This estimation shows that underflow may occur for computing the large zeros. In 
this case the range of computation of the weights can be enlarged by scaling the
factor  $e^{x^2/2}$ in the asymptotic approximations and computing scaled weights given by

\begin{equation}
\label{wghsHerm}
\tilde{w}_i =w_i e^{x^2_i}\,.
\end{equation}
With this, the overflow/underflow limitations are eliminated.

Using (\ref{wgh}) this scaled weight can be written as
\begin{equation}
\label{scaledw}
\tilde{w}_i=\Frac{\sqrt{\pi}2^{n+1}n!}{y'(x_i)^2},\quad y(x)=e^{-x^2/2}H_n(x)\,.
\end{equation}
This expression does not have overflow/underflow limitations neither with respect to $x$
nor with respect to $n$. Using (\ref{eq:herelem03}) or (\ref{eq:HermAiry01}) we observe that
the dominant factors $e^{-x^2/2}$ and $2^{n+1}n!$ can be explicitly cancelled out. 

Another
interesting property of this expression is that it is well conditioned with respect to the values
of the nodes. Indeed, we have $\tilde{w}_i=W(x_i)$, where we define the function 
$W(x)=\sqrt{\pi}2^{n+1}n!/y'(x)^2$. Now, it is straightforward to check that $W'(x_i)=0$
which means that, at the nodes $x=x_i$, the value of the weight is little affected by variations
on the actual value of the node. This, as we will show, will allow us to compute scaled weights
with nearly full double precision in all the range. 

For computing the scaled weights in this way, we need to compute $y'(x)$ from the asymptotic expansions
(\ref{eq:herelem03}) or (\ref{eq:HermAiry01}). This is a straightforward computation and, for instance,
starting from  (\ref{eq:herelem03}) we have that
\begin{equation}
y'(x)=\frac{2^{\frac12 n+1}g (\mu)}{\mu\left(1-t^2\right)^{5/4}}\left[
\cos\left(\mu^2 \eta-\tfrac14\pi\right) {\cal C}_{\mu}(t)-\sin\left(\mu^2 \eta-\tfrac14\pi\right) {\cal D}_{\mu}(t)\right],
\end{equation}
where
\begin{equation}
\begin{array}{l}
{\cal C}_{\mu}(t)\sim \displaystyle\sum_{s=0}^{\infty}\Frac{(-1)^s a_s}{(1-t^2)^{3s}\mu^{4s}},\quad
{\cal D}_{\mu}(t)\sim \displaystyle\sum_{s=-1}^{\infty}\Frac{(-1)^s b_s}{(1-t^2)^{3s+3/2}\mu^{4s+2}}
\end{array}
\end{equation}
and
\begin{equation}
\begin{array}{l}
a_s=\left(\frac12 +6s\right) t u_{2s} +u_{2s+1}+(1-t^2)\dot{u}_{2s},\\
\\
b_{-1}=1,b_s=\left(\frac72 +6s\right) t u_{2s+1} +u_{2s+2}+(1-t^2)\dot{u}_{2s+1},\,s\ge 0.
\end{array}
\end{equation}
The dots mean derivative with respect to $t$.

\begin{figure}
\begin{center}
\hspace*{-2cm}
\epsfxsize=16cm \epsfbox{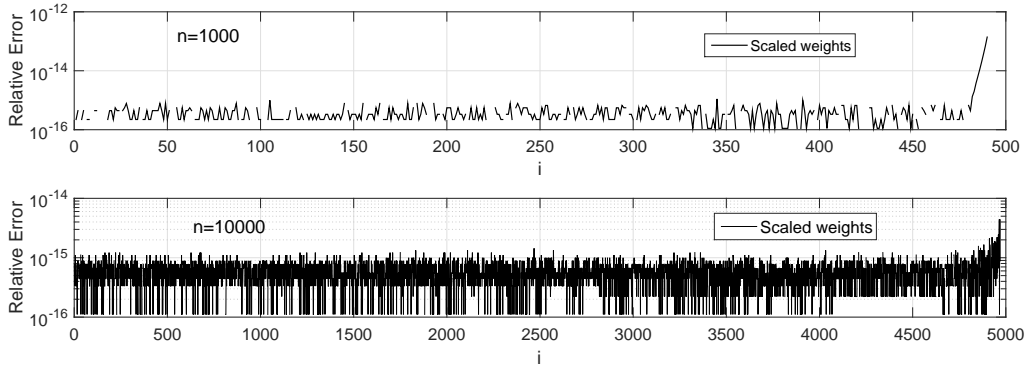}
\caption{ Relative accuracy obtained for the computation of the scaled weights (\ref{scaledw})
for $n=1000,\,10000$ using the asymptotic expansion for the Hermite polynomials
 in terms of elementary functions. 
\label{FigweightsH}}
\end{center}
\end{figure}

Two examples of computation of the scaled weights (for $n=1000,\,10000$) using
the expansion in terms of elementary functions are shown in Figure~\ref{FigweightsH}. 
As can be seen, most of the scaled weights can be computed with almost double precision accuracy. 
Also, as expected, there is some loss of accuracy for the weights corresponding to the largest nodes
(as discussed, for these values one has to use the expansion for the Hermite
polynomials in terms of Airy functions).

Typically, the
 additional computation of the weights requires about $70\%$ more CPU time than when using 
the asymptotic expansion in terms of elementary functions and about $133\%$ more CPU time
than when using the asymptotic expansion in terms of Airy functions (due to the computation of these
functions). This shows that, when possible, the direct computation of nodes and weights 
using asymptotics will be more efficient than computing more crude first approximations
and then refining with an iterative method which uses values of the orthogonal
polynomial. Each time the function (and its derivative
when we use Newton's method) is computed, the CPU time increases by this same amount, 
and only when one iteration is needed the speed would be comparable.

\section{Laguerre polynomials}\label{sec:Lagpol}

We consider asymptotic expansions for the Laguerre polynomials $L_n^{(\alpha)}(x)$ in terms of Bessel functions, Airy functions and Hermite polynomials.
Some of these expansions have been used to build an efficient scheme for computing the Laguerre polynomials for large values of $n$ and
small values of $\alpha$ ($-1<\alpha\le 5)$ \cite{gil:2016:ECL}. We discuss how to
use the expansions to obtain approximations to the zeros of Laguerre polynomials. 
Later, in Section \ref{lalp} we give expansions valid for large $n$ and $\alpha$.

For a survey of the work of several authors on inequalities and asymptotic formulas for the zeros of $L_n^{(\alpha)}(x)$ as $n$ or $\alpha$  or $\nu=4n+2\alpha+2$ $\to\infty$, we refer to \cite{Gatteschi:2002:ABL}.
 See also \cite{Huybrechs:2017:CAI}, were an alternative method, based on nonlinear steepest descent 
analysis of Riemann--Hilbert problems, is given for Laguerre-type Gaussian quadrature (and in
particular Gauss--Laguerre).

\subsection{A simple Bessel-type expansion}\label{sec:Lagsim}
We have the following representation\footnote{We summarize the results of  \cite[\S10.3.4]{Temme:2015:AMI}.}
\begin{equation}\label{eq:lagsim01}
L_n^{(\alpha)}(x)=\left(\frac{x}{n}\right)^{-\frac12\alpha} e^{\frac12x}
\left(J_{\alpha}\left(2\sqrt{nx}\right)A(x)-\sqrt{\frac{x}{n}}J_{\alpha+1}\left(2\sqrt{nx}\right)B(x)\right),
\end{equation}
with expansions
\begin{equation}\label{eq:lagsim02}
A(x)\sim\sum_{k=0}^\infty (-1)^k\frac{a_k(x)}{n^k},\quad B(x)=\sum_{k=0}^\infty (-1)^k \frac{b_k(x)}{n^k}\quad n\to\infty,
\end{equation}
valid for bounded values of $x$ and $\alpha$.

The coefficients $a_k(x)$ and $b_k(x)$ follow from the expansion of the function
\begin{equation}\label{eq:lagsim03}
f(z,s)=e^{xg(s)}\left(\frac{s}{1-e^{-s}}\right)^{\alpha+1},\quad g(s)=\frac{1}{s}-\frac{1}{e^s-1}-\frac12.
\end{equation}
The function $f$ is analytic in the strip $\vert\Im s\vert<2\pi$ and it can be expanded for  $\vert s\vert<2\pi$ into
\begin{equation}\label{eq:lagsim04}
f(x,s)=\sum_{k=0}^\infty c_k(x) s^k.
\end{equation}
The coefficients $c_k(x)$  are combinations of Bernoulli numbers and Bernoulli polynomials, the first ones being (with $c=\alpha+1$)
\begin{equation}\label{eq:lagsim05}
\begin{array}{@{}r@{\;}c@{\;}l@{}}
c_0(x)&=&1,\quad c_1(x)=\frac{1}{12}\left(6c-x\right),\quad \\[8pt]
c_2(x)&=&\frac{1}{288}\left(-12c+36c^2-12xc+x^2\right),\\[8pt]
c_3(x)&=&\frac{1}{51840}\left(-5x^3 + 90x^2c +(-540c^2 + 
 180c+72)x +1080c^2(c-1)\right).
\end{array}
\end{equation}
The coefficients $a_k(x)$ and $b_k(x)$ are in terms of the $c_k(x)$ given by
\begin{equation}\label{eq:lagsim06}
\begin{array}{@{}r@{\;}c@{\;}l@{}}
a_k(x) & = & \dsp{\sum_{m=0}^k \binomial{k}{m}(m+1-c)_{k-m}x^m c_{k+m}(x),}\\[8pt]
b_k(x) & = & \dsp{\sum_{m=0}^k \binomial{k}{m}(m+2-c)_{k-m}x^m c_{k+m+1}(x),}
\end{array}
\end{equation}
$k=0,1,2,\ldots$, and the first relations are
\begin{equation}\label{eq:lagsim07}
\begin{array}{ll}
a_0(x)= c_0(x)=1,\quad b_0(x)= c_1(x),\\[8pt]
a_1(x)= (1-c)c_1(x)+xc_2(x),\quad b_1(x)= (2-c)c_2(x)+xc_3(x),\\[8pt]
a_2(x)= (c^2-3c+2)c_2(x)+(4x-2xc)c_3(x)+x^2c_4(x),\\[8pt]
b_2(x)=  (c^2-5c+6)c_3(x)+(6x-2xc)c_4(x)+x^2c_5(x),
\end{array}
\end{equation}
again with $c=\alpha+1$.

\subsubsection{Expansions of the zeros}\label{sec:Lagsimzer}

Approximations of the zeros of $L_n^{(\alpha)}(x)$ can be obtained from \eqref{eq:lagsim01} and expressed in
terms of zeros of the Bessel function $J_\alpha(x)$. Because the expansion is valid for bounded values of $x$, 
the approximation can only be used for the small zeros. For example, in Table~\ref{tab:tab01} 
we show the results for the first 10 zeros when $n=100$, and for these early zeros the approximations are satisfactory.

We write (see \eqref{eq:lagsim01})
\begin{equation}\label{eq:lagsim08}
W(x)=J_{\alpha}\left(2\sqrt{nx}\right)A(x)-\sqrt{\frac{x}{n}}J_{\alpha+1}\left(2\sqrt{nx}\right)B(x),
\end{equation}
A first approximation to the zero $x_k$ of  $L_n^{(\alpha)}(x)$ follows from writing $2\sqrt{nx_k}=j_k$, where $j_k$ is the $k$th zero of $J_\alpha(x)$. A further approximation will be obtained by writing 
\begin{equation}\label{eq:lagsim09}
x_k=\xi+\eps, \quad \xi=\frac{1}{4n}j_k^2.
\end{equation}
By  expanding $W(x)$ at the zero $x=\xi+\eps$, assuming that $\eps$ is small, we find 
\begin{equation}\label{eq:lagsim10}
W(\xi)+\frac{\eps}{1!}W^\prime(\xi)+\frac{\eps^2}{2!}W^{\prime\prime}(\xi)+\ldots=0,
\end{equation}
and substituting an expansion of the form
\begin{equation}\label{eq:lagsim11}
\eps\sim\frac{\xi_1}{n}+\frac{\xi_2}{n^2}+\frac{\xi_3}{n^3}+\ldots,
\end{equation}
we find the following first few values
\begin{equation}\label{eq:lagsim12}
\begin{array}{@{}r@{\;}c@{\;}l@{}}
\xi_1&=&\dsp{ \frac{\xi}{12}}(\xi-6(\alpha+1)),\\[8pt]
\xi_2&=& 
\dsp{\frac{\xi}{720}}(150-90\xi+11\xi^2+360\alpha+210\alpha^2-90\xi\alpha),\\[8pt]
\xi_3&=& 
\dsp{\frac{\xi}{20160}}(2121\xi-770\xi^2+73\xi^3-6300\alpha-8820\alpha^2\,+\\[8pt]
&&\quad\quad 5040\xi\alpha-3780\alpha^3-770\xi^2\alpha+2751\xi\alpha^2-1260),
\end{array}
\end{equation}
where $\xi$ is defined in \eqref{eq:lagsim09}.

\paragraph{Algorithm and first numerical examples for the zeros}

The algorithm for computing the asymptotic approximation of the zeros runs in the same way as described for the Hermite polynomials, but is quite simple now. First compute $\xi$ from \eqref{eq:lagsim09} and the $\xi_j$ given in \eqref{eq:lagsim12}, then compute $\eps$ from \eqref{eq:lagsim11}, and finally $x_k$ from \eqref{eq:lagsim09}.

In Table~\ref{tab:tab01} we show the results of a first numerical verification for the expansion.
We take $n=100$, $\alpha=\frac13$, and compute the first 10 zeros by using Maple with Digits = 32. 
We show the relative errors in our approximations when we take 2, 4 and 6 terms in the expansion \eqref{eq:lagsim11}.
As can be seen in the table,  it is possible to obtain an accuracy near double precision ($10^{-16}$) in the computation of the first two zeros
of $L^{(1/3)}_{100}(x)$ using just the expansion with 6 terms. 

\renewcommand{\arraystretch}{1.2}
\begin{table}
\caption{
Relative errors in the computation of the zeros $x_k$ (see \eqref{eq:lagsim09}) by using the expansion \eqref{eq:lagsim11} with 2, 4 and 6 terms.
We take $n=100$, $\alpha=\frac13$.
\label{tab:tab01}}
$$
\begin{array}{rrrr}
k & 2\  {\rm terms} \ & 4\   {\rm terms} \quad &  6\  {\rm terms} \quad \\
\hline
1  &   0.22\times 10^{-6}  &   0.42\times 10^{-11}  &   0.21\times 10^{-15} \\
2  &   0.20\times 10^{-6}  &   0.34\times 10^{-11}  &   0.13\times 10^{-15} \\
3  &   0.18 \times 10^{-6}  &   0.22 \times 10^{-11}  &   0.84 \times 10^{-15} \\
4  &   0.15 \times 10^{-6}  &   0.96 \times 10^{-12}  &   0.22 \times 10^{-14}\\
5  &   0.12 \times 10^{-6}  &   0.11\times 10^{-13}  &   0.10 \times 10^{-14}\\
6  &   0.82\times 10^{-7}  &  0.43 \times 10^{-12}  &   0.67\times 10^{-14}\\
7  &   0.47\times 10^{-7}  &  0.29 \times 10^{-12}  &   0.24\times 10^{-13}\\
8  &   0.16\times 10^{-7}  &   0.28\times 10^{-12}  &   0.48\times 10^{-13}\\
9  &   0.87\times 10^{-8}  &   0.95\times 10^{-12}  &   0.72\times 10^{-13}\\
10  &   0.24\times 10^{-7}  &   0.13\times 10^{-11}  &   0.87\times 10^{-13}\\
\hline
\end{array}
$$
\end{table}
\renewcommand{\arraystretch}{1.0}

\subsection{An expansion in terms of Airy functions}\label{sec:LagAiryexp}
We start with the representation\footnote{We summarize results of \cite{Frenzen:1988:UAE};  see also \cite[\S{VII.5}]{Wong:2001:AAI}.}
\begin{equation}\label{eq:lagairy01}
L_n^{(\alpha)}(\nu \sigma)=(-1)^n
\frac{ e^{\frac12\nu \sigma}\chi(\zeta)}{2^\alpha\nu^{\frac13}}\left(\Ai\left(\nu^{2/3} \zeta\right)A(\zeta)
+\nu^{-\frac43}\Ai^{\prime}\left(\nu^{2/3}\zeta\right)
B(\zeta)\right)
\end{equation}
with  expansions
\begin{equation}\label{eq:lagairy02}
A(\zeta)\sim\sum_{j=0}^\infty\frac{\alpha_{2j}}{\nu^{2j}},\quad B(\zeta)\sim\sum_{j=0}^\infty\frac{\beta_{2j+1}}{\nu^{2j}},\quad n\to\infty,
\end{equation}
uniformly for bounded $\alpha$ and  $\sigma\in(\sigma_0,\infty]$, where $\sigma_0\in(0,1)$, a fixed number. 

Here
\begin{equation}\label{eq:lagairy03}
\nu=4\kappa,\quad \kappa=n+\tfrac12(\alpha+1), \quad 
\chi(\zeta)=2^{\frac12}\sigma^{-\frac14-\frac12\alpha}\left(\frac{\zeta}{\sigma-1}\right)^{\frac14},
\end{equation}
and
\begin{equation}\label{eq:lagairy04}
\begin{cases}
\dsp{\tfrac23(-\zeta)^{\frac32}=\tfrac12\left(\arccos\sqrt{{\sigma}}-\sqrt{{\sigma-\sigma^2}}
\right)}& \quad\text{if $\quad 0<\sigma\le1$,}
\\[8pt]
\dsp{\tfrac23\zeta^{\frac32}=\tfrac12\left(\sqrt{{\sigma^2-\sigma}}-\arccosh\sqrt{{\sigma}}\right)}& 
\quad\text{if $\quad \sigma\ge1$.}
\end{cases}
\end{equation}
We have the relation
\begin{equation}\label{eq:lagairy05}
\zeta^\frac12\frac{d\zeta}{d\sigma}=\frac{\sqrt{\sigma-1}}{2\sqrt{\sigma}}.
\end{equation}

For the derivative we can use the relation
\begin{equation}\label{eq:lagairy06}
\frac{d}{dx}L_n^{(\alpha)}(x)=L_n^{(\alpha)}(x)-L_n^{(\alpha+1)}(x).
\end{equation}

The first coefficients of the expansions in \eqref{eq:lagairy02} are
\begin{equation}\label{eq:lagairy07}
\alpha_0=1,\quad\ \beta_1=-\frac{1}{4b^3}\left(f_1-bf_2\right),
\end{equation}
where $b=\sqrt{\zeta}$ if $\zeta\ge0$ and $b=i\sqrt{-\zeta}$ when $\zeta\le0$, and
\begin{equation}\label{eq:lagairy08}
\begin{array}{ll}
\dsp{f_1=i\frac{\left(\sigma+3\alpha(\sigma-1)\right)\sigma^2a_1^3-2}{3a_1^2\sigma\sqrt{\sigma(1-\sigma)}}
,}\\[8pt] 
\dsp{f_2=\frac{-4-8\sigma^2(\sigma+3\sigma\alpha-3\alpha)a_1^3+\sigma^4(12\sigma-3-4\sigma^2+12\alpha^2(\sigma-1)^2)a_1^6}{12\sigma^3a_1^4(\sigma-1)},}\\[8pt] 
\dsp{a_1=\left(\frac{4\zeta}{\sigma^3(\sigma-1)}\right)^{\frac14}}.
\end{array}
\end{equation}
More coefficients can be obtained by the method described in \cite[\S23.2]{Temme:2015:AMI}. Starting point in this case is the integral (see 
\cite[\S{VII.5, (5.11)}]{Wong:2001:AAI})
\begin{equation}\label{eq:lagairy09}
\frac{1}{2\pi i}\int_{\calL} f(u) e^{\nu\left(\frac13u^3-\zeta u\right)}\,du,
\end{equation}
where $\calL$ is an Airy-type contour and $f(u)$  is given by 
\begin{equation}\label{eq:lagairy10}
f(u)=(1-z^2)^{\frac12(\alpha-1)}\frac{dz}{du}.
\end{equation}
The relation between $z$ and $u$ follows in this case from the cubical transformation
\begin{equation}\label{eq:lagairy11}
\tfrac12\arctanh \,z-\tfrac12 z\sigma
=\tfrac13u^3-\zeta u, \quad \frac{dz}{du}= \frac{2(u^2-\zeta)(1-z^2)}{1-\sigma+\sigma z^2}.
\end{equation}
The function $f(u)$ can be expanded in a two-point Taylor series 
\begin{equation}\label{eq:lagairy12}
f(u)=\sum_{k=0}^\infty(c_k+ud_k)\left(u^2-\zeta\right)^k,
\end{equation}
in which the coefficients can be expressed in terms of the derivatives of $f(u)$ at $u=\pm\sqrt{\zeta}$. An integration by parts procedure then gives the coefficients $\alpha_{2j}$ and $\beta_{2j+1}$ of \eqref{eq:lagairy02}.

In \S\ref{sec:Bessalgcoef} we describe in detail this method for a Bessel-type expansion.

\subsubsection{Expansions of the zeros}\label{sec:LagAiryzer}
We write
\begin{equation}\label{eq:lagairy13}
W(\zeta)=\Ai\left(\nu^{2/3} \zeta\right)A(\zeta)+\nu^{-\frac43}\Ai^{\prime}\left(\nu^{2/3}\zeta\right)B(\zeta),
\end{equation}
where $A(\zeta)$ and $B(\zeta)$  have the expansions shown in \eqref{eq:lagairy02}. 

Similarly as in \S\ref{sec:Hermpol} we write the zeros $x_j$ of 
$L_n^{(\alpha)}(x)$ in terms of the zeros $a_k$ of the Airy function. These zeros are negative, and $a_1$ will correspond the $n$th zero of $L_n^{(\alpha)}(x)$, $a_2$ with the $(n-1)$th zero, and so on. 

A zero of  $L_n^{(\alpha)}(x)$ is a zero of $W(\zeta)$ and it can be written in terms of $\zeta$ in the form
\begin{equation}\label{eq:lagairy14}
\zeta=\zeta_0+\eps, \quad \zeta_0=\nu^{-\frac23}a_j, \
\end{equation}
and we assume that we can expand
\begin{equation}\label{eq:lagairy15}
\eps\sim\frac{\zeta_1}{\nu^2}+\frac{\zeta_2}{\nu^4}+\frac{\zeta_3}{\nu^6}+\ldots.
\end{equation}
By  expanding $W(\zeta)$ at $\zeta_0$ we have
\begin{equation}\label{eq:lagairy16}
W(\zeta_0)+\frac{\eps}{1!}W^\prime(\zeta_0)+\frac{\eps^2}{2!}W^{\prime\prime}(\zeta_0)+\ldots=0,
\end{equation}
and substituting the expansions shown in  \eqref{eq:lagairy02} we can obtain the coefficients $\zeta_j$. We obtain
\begin{equation}\label{eq:lagairy17}
\zeta_1=-\beta_1,\quad \zeta_2=-\left(\beta_3+\tfrac16\zeta_0\zeta_1^3+\zeta_1\alpha_2+\zeta_1\frac{d}{d\zeta}\beta_1+\tfrac12\zeta_0\beta_3\zeta_1^2\right),
\end{equation}
where $\beta_1$ given in \eqref{eq:lagairy07}. The coefficients are evaluated at $\zeta_0$.

\paragraph{Algorithm and first numerical examples for the zeros}
In \S\ref{sec:Hermsimzer} we have described the algorithm for computing the asymptotic approximation of the zeros for the Airy case. The present algorithm runs in the same way. 
For the zero $x_{n+1-j}$, $j=1,2,\ldots$, first compute $\zeta_0$, from \eqref{eq:lagairy14}. Then 
compute $\sigma_0$ by inverting the first relation in \eqref{eq:lagairy04}. This is done by using the expansion
\begin{equation}\label{eq:lagairy18}
\sigma=1+\wt\zeta+\tfrac15\wt\zeta^2-\tfrac{3}{175}\wt\zeta^3+\tfrac{23}{7875}\wt\zeta^4+\ldots,\quad \wt\zeta=2^{\frac23}\zeta.\end{equation}
An alternative would be to use an iterative method. In that case it is convenient to write $\sigma=\cos^2\theta$, and the equation to be solved for $\theta$ becomes $\frac83(-\zeta)^{\frac32}=\theta-\sin\theta$, $0\le\theta<\pi$. 

With  $\sigma=\sigma_0$ we compute the coefficients in \eqref{eq:lagairy17}, then $\eps$ and $\zeta$ from \eqref{eq:lagairy15} and \eqref{eq:lagairy14}. A final inversion of the relation  in the first line of \eqref{eq:lagairy04} gives the $\sigma$, and then $x_{n+1-j}\sim \nu \sigma$.

For example, we take $n=100$, $\alpha=\frac13$, and we compute the zero $x_{100}=375.635158667\cdots$ by using Maple. We compute (see  \eqref{eq:lagairy14})
\begin{equation}\label{eq:lagairy19}
\zeta_0= -2.3381074105\nu^{-2/3}=-0.0428779491924.
\end{equation}
Upon solving the first equation in \eqref{eq:lagairy04} for $\sigma$, we obtain  
$\sigma_0=0.9328675228515$. With this value, a first approximation of the zero is  
$x_{100}\sim \nu \sigma_0\doteq 375.634655868$, with a relative accuracy  of $1.34 \times 10^{-6}$.

Finally, we compute $\zeta_1=0.131145197575$, compute $\zeta\sim \zeta+\zeta_1/\nu^2$, invert again the first relation in \eqref{eq:lagairy04}, giving
$\sigma=0.932868771534$ and $x_{100}\doteq 375.635158671$, a relative accuracy  of $1.08 \times 10^{-11}$.
 For the halfway zero $x_{51}$ we found the relative accuracies  $6.71 \times 10^{-6}$ and $1.52 \times 10^{-10}$.

\subsection{Another expansion in terms of Bessel functions}\label{sec:LagBesexp}
After substituting $t=e^{-s}$ in the integral representation\footnote{We summarize the results of \cite{Frenzen:1988:UAE}; see also \cite[\S{VII.7}]{Wong:2001:AAI}.} 
\begin{equation}\label{eq:lagbess01}
L_n^{(\alpha)}(z)= \frac{1}{2\pi i}\int_\calL (1-t)^{-\alpha-1}e^{-tz/(1-t)}\,\frac{dt}{t^{n+1}},
\end{equation}
we obtain the representation
\begin{equation}\label{eq:lagbess02}
e^{-\nu \rho}L_n^{(\alpha)}(2\nu \rho)=\frac{2^{-\alpha}}{2\pi i}\int_{-\infty}^{(0+)}e^{\nu h(s,\rho)}\left(\frac{\sinh s}{s}\right)^{-\alpha-1}\,\frac{ds}{s^{\alpha+1}},
\end{equation}
where $\nu=2n+\alpha+1$ and $h(s,\rho)=s-\rho\coth s$. The contour  starts at $-\infty$ with  $\ph\,u=-\pi$, encircles the origin anti-clockwise, and returns to $-\infty$ with $\ph\, u=\pi$.  The transformation to a standard form for this case is $h(s,\rho)=u-\zeta/u$, with result
\begin{equation}\label{eq:lagbess03}
2^{\alpha}e^{-\nu \rho}L_n^{(\alpha)}(2\nu \rho)=\frac{1}{2\pi i}\int_{-\infty}^{(0+)}e^{\nu(u-\zeta/u)}f(u)\,\frac{du}{u^{\alpha+1}},
\end{equation}
where
\begin{equation}\label{eq:lagbess04}
f(u)=\left(\frac{u}{\sinh s}\right)^{\alpha+1}\frac{ds}{du}.
\end{equation}
By using an integration by parts procedure (see \S\ref{sec:Bessalgcoef}), we can obtain the representation
\begin{equation}\label{eq:lagbess05}
L_n^{(\alpha)}(2\nu \rho)= \frac{e^{\nu \rho}\chi(\zeta)}{2^\alpha\zeta^{\frac12\alpha}}\left(J_\alpha\bigl(2\nu \sqrt{\zeta}\bigr)
A(\zeta)-
\frac{1}{\sqrt{\zeta}}J_{\alpha+1}\bigl(2\nu  \sqrt{\zeta}\bigr)B(\zeta)\right),
\end{equation}
with expansions
\begin{equation}\label{eq:lagbess06}
A(\zeta)\sim\sum_{j=0}^\infty\frac{A_{2j}(\zeta)}{\nu^{2j}},\quad 
B(\zeta)\sim\sum_{j=0}^\infty\frac{B_{2j+1}(\zeta)}{\nu^{2j+1}},\quad \nu\to\infty,
\end{equation}
uniformly for $\rho\le 1-\delta$, where $\delta\in(0,1)$  is a fixed number. Here,
\begin{equation}\label{eq:lagbess07}
\nu=2n+\alpha+1,\quad \chi(\zeta)=(1-\rho)^{-\frac14}\left(\frac{\zeta}{\rho}\right)^{\frac12\alpha+\frac14},\quad \rho<1,
\end{equation}
with $\zeta$ given by
\begin{equation}\label{eq:lagbess08}
\begin{cases}
\dsp{\sqrt{-\zeta}=\tfrac12\left(\sqrt{{\rho^2-\rho}}
+\arcsinh\sqrt{{-\rho}}\right)}, & \quad\text{if \quad$\rho\le0$,}
\\[8pt]
\dsp{\sqrt{\zeta}=\tfrac12\left(\sqrt{{\rho-\rho^2}}+\arcsin\sqrt{{\rho}}\right)}, 
& \quad \text{if \quad$ 0\le \rho<1$.}
\end{cases}
\end{equation}
 We have the relation
\begin{equation}\label{eq:lagbess09}
\frac{1}{\zeta^{\frac12}}\frac{d\zeta}{d\rho}=\sqrt{\frac{1-\rho}{\rho}},\quad \rho <1.
\end{equation}

The first coefficients are
\begin{equation}\label{eq:lagbess10}
\begin{array}{@{}r@{\;}c@{\;}l@{}}
A_0(\zeta)&=&1, \\[8pt] 
B_1(\zeta)&=&\dsp{ \frac{1}{48\xi}\left(5\xi^4b+6\xi^2b+3\xi+12a^2(b-\xi)-3b\right),}
\end{array}
\end{equation}
where
\begin{equation}\label{eq:lagbess11}
\xi=\sqrt{\frac{\rho}{1-\rho}}, \quad b=\sqrt{\zeta}.
\end{equation}
More coefficients can be obtained by using the method described in  \S\ref{sec:Bessalgcoef}.

To remove in \eqref{eq:lagbess05} the singularities due  to the Bessel functions at $\zeta=0$,  it is convenient to use the function
$E_\nu(z)$ introduced by Tricomi; see  \cite[p.~34]{Tricomi:1947:SFI}. We have
\begin{equation}\label{eq:lagbess12}
E_\nu(z)=z^{-\frac12\nu} J_\nu\left(2\sqrt{z}\right)=\sum_{k=0}^\infty(-1)^k\frac{z^k}{k!\,\Gamma(\nu+k+1)}.
\end{equation}
It is an analytic function of $z$. In terms of the modified Bessel  function we can write
\begin{equation}\label{eq:lagbess13}
E_\nu(-z)=z^{-\frac12\nu} I_\nu\left(2\sqrt{z}\right)=\sum_{k=0}^\infty\frac{z^k}{k!\,\Gamma(\nu+k+1)}.
\end{equation}
The representation in \eqref{eq:lagbess05} can be written in the form
\begin{equation}\label{eq:lagbess14}
L_n^{(\alpha)}(2\nu \rho)=\left(\tfrac12 \nu\right)^{\alpha} e^{\nu \rho}
\chi(\zeta)\Bigl(E_\alpha\bigl(\zeta\nu^2\bigr)
A(\zeta)-E_{\alpha+1}\bigl(\zeta\nu^2\bigr)B(\zeta)\Bigr),
\end{equation}
and we can use this representation also for $\zeta<0$, i.e., $\rho<0$.

For more details about the coefficients $A_j(\zeta)$ and $B_j(\zeta)$ of the expansions in \eqref{eq:lagbess06}, see  
\cite{gil:2016:ECL}.

\subsubsection{A general method for  the coefficients in Bessel-type expansions}\label{sec:Bessalgcoef}
We describe a general method for evaluating the coefficients  $A_k(\zeta)$ and $B_k(\zeta)$ used in \eqref{eq:lagbess06}. 

We consider the standard form
\begin{equation}\label{eq:Bessalgcoef01}
F_\zeta(\nu)=\frac1{2\pi i}\int_{{\calC}}e^{\nu\left(u-\zeta/u\right)}f(u)\,\frac{du}{u^{\alpha+1}},
\end{equation}
where the contour ${\calC}$ starts at $-\infty$ with  $\ph\,u=-\pi$, encircles the origin anti-clockwise, and returns to $-\infty$ with $\ph\, u=\pi$. The $f(u)$ is assumed to be analytic in a neighborhood of $\calC$, and  in particular in a domain that contains the saddle points $\pm ib$, where $b=\sqrt{\zeta}$. 

When we replace $f$ by unity, we obtain the Bessel function:
\begin{equation}\label{eq:Bessalgcoef02}
F_\zeta(\nu)=\zeta^{-\frac12\alpha} J_\alpha\left(2\nu\sqrt{\zeta}\right). 
\end{equation}

The coefficients of the expansions in \eqref{eq:lagbess06} follow from the recursive scheme
\begin{equation}\label{eq:Bessalgcoef03}
\begin{array}{@{}r@{\;}c@{\;}l@{}}
f_j(u)&=&\dsp{A_j(\zeta)+B_j(\zeta)/u+\left(1+b^2/u^2\right)g_j(u),}\\[8pt]
f_{j+1}(u)&=&\dsp{g_j^\prime(u)-\frac{\alpha+1}{u}g_j(u),}\\[8pt]
A_j(\zeta)&=&\dsp{\frac{f_j(ib)+f_j(-ib)}{2}},\quad \dsp{B_j(\zeta)=i\frac{f_j(ib)-f_j(-ib)}{2b},}
\end{array}
\end{equation}
with $f_0(u)=f(u)$, the coefficient function. 

Using  this scheme and integration by parts, we can obtain the asymptotic expansion
\begin{equation}\label{eq:Bessalgcoef04}
F_\eta(\nu)\sim \zeta^{-\frac12\alpha} J_\alpha\left(2\nu\sqrt{\zeta}\right)\sum_{j=0}^\infty(-1)^j\frac{A_j(\zeta)}{\nu^j}+
 \zeta^{-\frac12(\alpha+1)} J_{\alpha+1}\left(2\nu\sqrt{\zeta}\right)\sum_{j=0}^\infty(-1)^j\frac{B_j(\zeta)}{\nu^j}.
 \end{equation}

The coefficients $A_j(\zeta)$ and $B_j(\zeta)$ can all be expressed in terms of the derivatives $f^{(k)}(\pm ib)$
 of $f(u)$ at the saddle points $\pm ib$; we will need these for $0\le k \le 2j$ (see \eqref{eq:Bessalgcoef08}).

We expand the functions $f_j(u)$ in two-point Taylor expansions
\begin{equation}\label{eq:Bessalgcoef05}
f_j(u)= \sum_{k=0}^\infty  C_k^{(j)} (u^2-b^2)^k+u\sum_{k=0}^\infty  D_k^{(j)} (u^2-b^2)^k.
\end{equation}

Using  \eqref{eq:lagbess03}, we derive the following recursive scheme for the coefficients
\begin{equation}\label{eq:Bessalgcoef06}
\begin{array}{ll}
\dsp{C_k^{(j+1)}=(2k-\alpha)D_{k}^{(j)}+b^2(\alpha-4k-2)D_{k+1}^{(j)}+2(k+1)b^4D_{k+2}^{(j)},}\\[8pt]
\dsp{D_k^{(j+1)}=(2k+1-\alpha)C_{k+1}^{(j)}-2(k+1)b^2C_{k+2}^{(j)},}
\end{array}
\end{equation}
for $j,k=0,1,2,\ldots$, and the coefficients $A_j$ and $B_j$ follow from 
\begin{equation}\label{eq:Bessalgcoef07}
A_j(\zeta)=C_0^{(j)},\quad B_j(\zeta)=-b^2D_0^{(j)},\quad j \ge 0.
\end{equation}

In the present case of the Laguerre polynomials the functions $f_{2j}$ are even and $f_{2j+1}$ are odd, and we have $A_{2j+1}(\zeta)=0$ and $B_{2j}(\zeta)=0$. A few non--vanishing coefficients are
\begin{equation}\label{eq:Bessalgcoef08}
\begin{array}{@{}r@{\;}c@{\;}l@{}}
A_0(\zeta)&=&f(ib),\\[8pt]
B_1(\zeta)&=&-\frac{1}{4}b\bigl((2\alpha-1)if^{(1)}(ib)+bf^{(2)}(ib)\bigr),\\[8pt]
A_2(\zeta)&=&-\frac{1}{32b}\bigl(3i(4^2\alpha-1)f^{(1)}(ib)-
(3-16\alpha+4\alpha^2)bf^{(2)}(ib)\ +\\[8pt]
&&2i(2\alpha-3)b^2f^{(3)}(ib)b^2+b^3f^{(4)}(ib)\bigr),\\[8pt]
B_3(\zeta)&=&-\frac{1}{384b}\bigl(3(4\alpha^2-1)(2\alpha-3)(if^{(1)}(ib)+bf^{(2)}(ib))\ +\\[8pt]
&&2i(\alpha-7)(2\alpha-1)(2\alpha-3)b^2f^{(3)}(ib)\ +\\[8pt]
&&3(19-20\alpha+4\alpha^2))b^3f^{(4)}(ib)
-3i(2\alpha-5)b^4f^{(5)}(ib)-b^5f^{(6)}(ib)\bigr).
\end{array}
\end{equation}

To have $A_0(\zeta)=1$ in  the first expansion in \eqref{eq:lagbess06} we have scaled all $A$ and $B$-coefficients with respect to $A_0(\zeta)=\chi(\zeta)$; see \eqref{eq:lagbess07}.

\begin{remark}\label{rem:twop}
The main step for obtaining the coefficients $A_j(\zeta)$ and $B_j(\zeta)$ is the evaluation of those for $j=0$ in \eqref{eq:Bessalgcoef05} and we summarize  the method described in \cite{Lopez:2002:TPT}. We rewrite the two-point Taylor expansion in the form
\begin{equation}\label{eq:Bessalgcoef09}
f(u)=\sum_{k=0}^\infty \bigl(
a_k(u_1,u_2)(u-u_1)+a_k(u_2,u_1)(u-u_2)\bigr)
(u-u_1)^k(u-u_2)^k,
\end{equation}
where, in the present case, $u_1=-b$ and $u_2=b$.
Then, 
\begin{equation}\label{eq:Bessalgcoef10}
C_k^{(0)}=-u_1a_k(u_1,u_2)-u_2a_k(u_2,u_1),\quad D_k^{(0)}= a_k(u_1,u_2)+a_k(u_2,u_1).
\end{equation}
We have $a_0(u_1,u_2)={f(b)}/{ (2b)}$\ and $ a_0(u_2,u_1)=-f(-b)/{ (2b)}$, and,  for $k=1,2,3,...$,
\begin{equation}\label{eq:Bessalgcoef11}
a_k(u_1,u_2)=
\sum_{j=0}^k\frac{(k+j-1)!}{
j!(k-j)!}\frac{(-1)^{k+1}kf^{(k-j)}(b)+(-1)^j jf^{(k-j)}(-b)}{
k!(-2b)^{k+j+1}},
\end{equation}
$a_k(u_2,u_1)$ follows from $a_k(u_1,u_2)$ by replacing $b$ by $-b$.
\eoremark
\end{remark}
\subsubsection{Expansions of the zeros}\label{sec:LagBeszer}
From the Bessel-type expansion we derive expansions of the first  half of the zeros of the Laguerre polynomial. We write
\begin{equation}\label{eq:LagBeszer01}
W(\zeta)=J_\alpha\bigl(2\nu \sqrt{\zeta}\bigr)A(\zeta)-\frac{1}{\sqrt{\zeta}}J_{\alpha+1}\bigl(2\nu  \sqrt{\zeta}\bigr)B(\zeta).
\end{equation}
A zero of  $L_n^{(\alpha)}(2\nu x)$ is a zero of $W(\zeta)$ and it can be written in terms of $\zeta$ in the form
\begin{equation}\label{eq:LagBeszer02}
\zeta=\zeta_0+\eps, \quad \zeta_0=\frac{j_k^2}{4\nu^2},
\end{equation}
where $j_k$ is a zero of $J_\alpha(z)$. By  expanding $W(\zeta)$ we have with the zero $\zeta$ in  this form 
\begin{equation}\label{eq:LagBeszer03}
W(\zeta_0)+\frac{\eps}{1!}W^\prime(\zeta_0)+\frac{\eps^2}{2!}W^{\prime\prime}(\zeta_0)+\ldots=0.
\end{equation}
We assume that $\eps$ can be expanded in the form
\begin{equation}\label{eq:LagBeszer04}
\eps\sim\frac{\zeta_1}{\nu^2}+\frac{\zeta_2}{\nu^4}+\frac{\zeta_3}{\nu^6}+\ldots,
\end{equation}
and substituting this expansion, we obtain $\zeta_1=-B_1(\zeta)$ (see \eqref{eq:lagbess10}) and
\begin{equation}\label{eq:LagBeszer05}
6\zeta\zeta_2=2 B_1(\zeta)^3-3 (\alpha+1) B_1(\zeta)^2+6 \zeta B_1(\zeta) \left(B_1^\prime(\zeta)+A_2(\zeta)\right)
-6 \zeta B_3(\zeta),
\end{equation}
In the algorithm we use $\zeta=\zeta_0$.

\paragraph{Algorithm and first numerical examples for the zeros}
As in the previous cases we describe how the asymptotic approximations for the zeros can be obtained. For the zero $x_{k}$, $k=1,2,\ldots$, first compute $\zeta_0$, from \eqref{eq:LagBeszer02}. Then compute 
$\rho_0$ by inverting the second relation in \eqref{eq:lagbess08}. This is done  by using the expansion
\begin{equation}\label{eq:LagBeszer06}
\rho=\zeta+\tfrac13\zeta^2+\tfrac{11}{45}\zeta^3+\tfrac{73}{315}\zeta^4+\ldots.
\end{equation}
An alternative is solving with an iterative method. In that case it is convenient to write 
$\rho=\sin^2\frac12\theta$, and the equation to be solved for $\theta$ becomes 
$8\sqrt{\zeta}=\theta+\sin\theta$, $0\le\theta<\pi$. With  
$\rho=\rho_0$ we compute the coefficients $\zeta_j$ in \eqref{eq:LagBeszer04}, see also \eqref{eq:lagbess10}. Compute $\zeta$ from \eqref{eq:LagBeszer02} and perform a final inversion of the relation  in the second line of \eqref{eq:lagbess08}. This gives the $\rho$, and then $x_{k}\sim 2\nu \rho$.

Because the expansions in \eqref{eq:lagbess06} become useless when $\rho\to1$, we should use the present result for a limited number of zeros, say, only for $k=1,2,3,\ldots,\frac12n$; the remaining zeros can be obtained by using the Airy-type expansion.

When we take $n=100$, $\alpha=\frac13$, and use the approximation $\zeta\sim \zeta_0$ with the first zero $x_{1}=0.02092331638663936$ computed by Maple with Digits=16, we found a relative accuracy of $3.65 \times 10^{-6}$; with the  term $\zeta_1/\nu^2$ included we found $6.68\times10^{-11}$
and when included up to the term $\zeta_3/\nu^6$, the accuracy is $2\times10^{-16}$. 
For the zero $x_{50}$ we found the relative errors $4.91 \times 10^{-6}$, $1.57 \times 10^{-10}$ and~$0$ (full double accuracy), respectively. 

In the next section we analyze in more detail the performance of the different expansions for the 
zeros and we also discuss the stable computation of the weights.

\subsection{Numerical performance of the expansions for $\alpha$ small}

\begin{figure}
\begin{center}
\hspace*{-2.5cm}
\epsfxsize=17cm \epsfbox{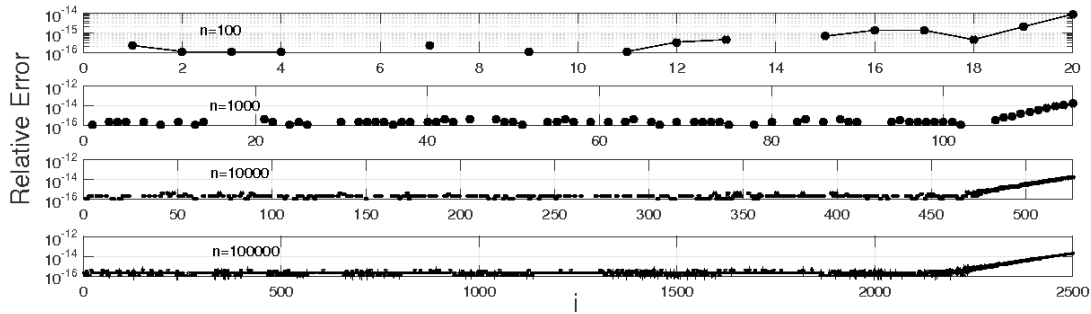}
\caption{ Relative accuracy obtained with the asymptotic expansion (\ref{eq:lagsim09}) for computing the first zeros of $L^{(1/4)}_n(x)$
for $n=100,\,1000,\,10000,\,100000$. 
The
label $i$ in the abscissa represents the order of the zero (starting from $i=1$ for the smallest zero). 
The points not shown in the plots correspond to values with all digits correct in double precision accuracy.
\label{Fig3}}
\end{center}
\end{figure}

\begin{figure}
\begin{center}
\hspace*{-2cm}
\epsfxsize=16cm \epsfbox{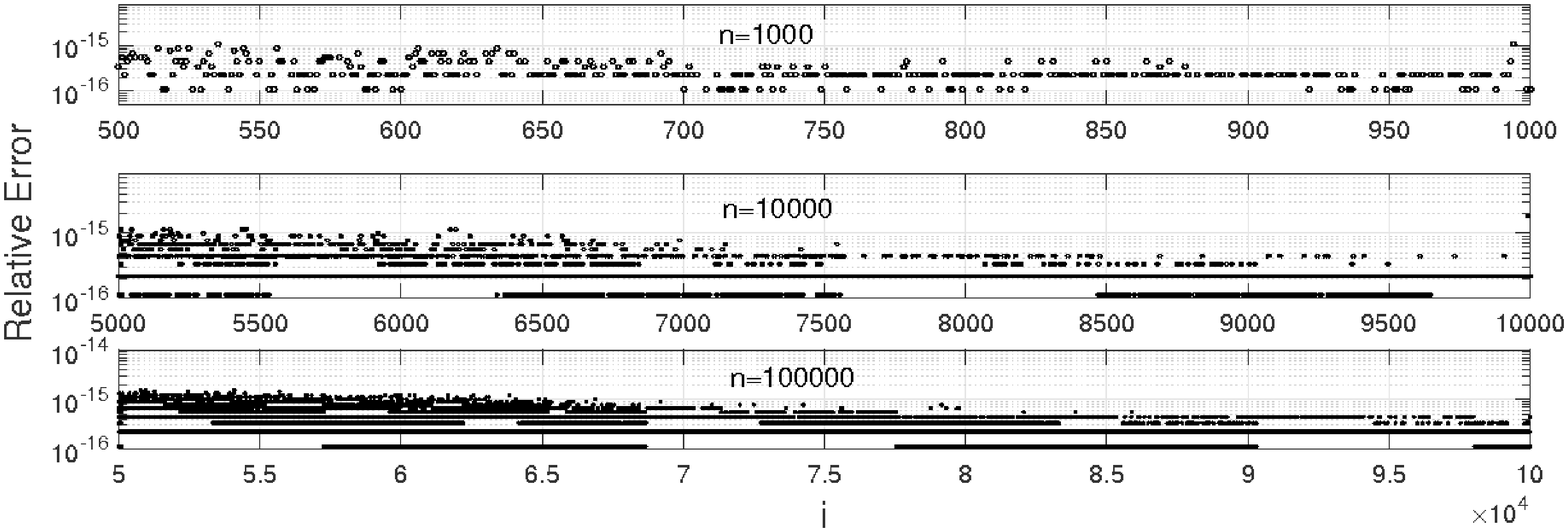}
\caption{ Relative accuracy obtained with the asymptotic expansion  (\ref{eq:lagairy18}) for computing the large zeros of $L^{(1/4)}_n(x)$
for $n=1000,\,10000,\,100000$.
\label{Fig4}}
\end{center}
\end{figure}

\begin{figure}
\begin{center}
\hspace*{-2cm}
\epsfxsize=16cm \epsfbox{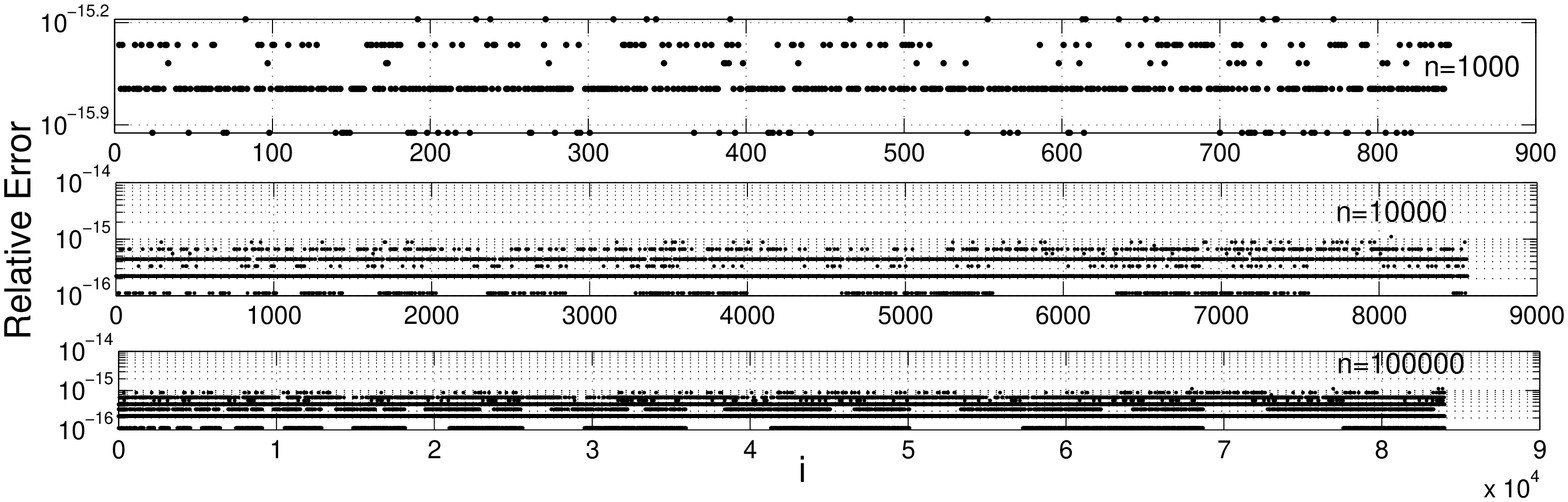}
\caption{ Relative accuracy obtained with the asymptotic expansion (\ref{eq:LagBeszer02}) for computing the zeros of $L^{(1/4)}_n(x)$
for $n=1000,\,10000,\,100000$.
\label{Fig5}}
\end{center}
\end{figure}

In Figures \ref{Fig3}, \ref{Fig4} and \ref{Fig5} we show the accuracy obtained with
the asymptotic expansions (\ref{eq:lagsim09}), (\ref{eq:lagairy18}), (\ref{eq:LagBeszer02}), respectively, for the zeros of the Laguerre polynomial $L^{(1/4)}_n(x)$ 
for different values of $n$. An implementation of the expansions in finite precision arithmetic (coded in Fortran 90) 
has been considered for testing. As for the Hermite case, in these implementations
only non-iterative methods (power series) are used for the inversion of the variables. For computing
the first zeros of Bessel functions we use the algorithm describe in \cite{Gil:2012:CZC}. For large zeros we use the MacMahon's
expansion (see \cite[\S10.21(vi)]{Olver:2010:Bessel})

\begin{equation}
j_{\nu,m} \sim a-\frac{\mu-1}{8a}-\frac{4(\mu-1)(7\mu-31)}{3(8a)^{3}}%
-\frac{32(\mu-1)(83\mu^{2}-982\mu+3779)}{15(8a)^{5}}-\cdots,.
\end{equation}
where $\mu=4\nu^2$, $a=(m+\nu/2-1/4)\pi$.
  
\begin{figure}
\begin{center}
\hspace*{-1.0cm}
\epsfxsize=14cm \epsfbox{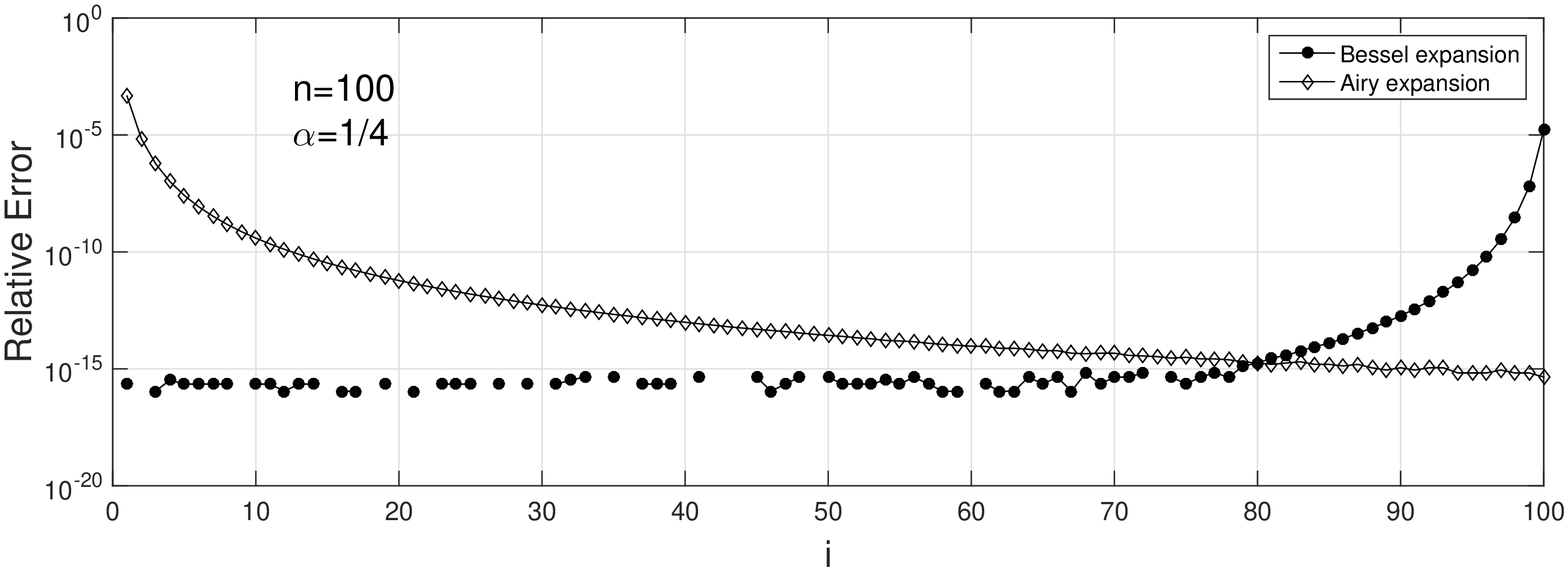}
\caption{ Relative accuracy obtained with the asymptotic expansions (\ref{eq:lagairy18}) and 
 (\ref{eq:LagBeszer02}) for computing the zeros of $L^{(1/4)}_n(x)$
for $n=100$. The points not shown in the plot correspond to values with all digits correct in double precision accuracy.
\label{Fig6}}
\end{center}
\end{figure}

As can be seen in Figure  \ref{Fig3}, the validity of the first asymptotic expansion 
in terms of zeros of Bessel functions (\ref{eq:lagsim11}) is limited to the first zeros.
 On the contrary, Figure \ref{Fig5}
shows that the other Bessel expansion (\ref{eq:LagBeszer02}) works very well for approximating a large number of zeros 
of the Laguerre polynomial but fails for the last zeros. For these zeros, the Airy expansion (\ref{eq:lagairy18})
should be used. The accuracy of the Bessel and Airy expansions for $n=100$ is illustrated in  Figure \ref{Fig6}.
As in the case of the Hermite approximations, the combined use of the expansions allow the computation of the
zeros of Laguerre polynomials for $n=100$ with an accuracy of $15$-$16$ digits.

The efficiency of the expansions is compared in Table \ref{tab:tab03}.
As in the case of the zeros of Hermite polynomials,
 in order to improve the speed of the methods we apply the expansions only 
in the regions where the inversion of the variables can be done accurately 
by using the series expansions 
(\ref{eq:lagairy18}) and 
(\ref{eq:LagBeszer06})
 in the case of the Airy expansion and the second Bessel expansion 
(\ref{eq:LagBeszer02}), respectively: the first $0.75n$ zeros for the Bessel expansion
and the last $0.25n$ zeros for the Airy expansion. 
For these two expansions, we observe in the Table
that there is no much difference in speed between using $2$ terms and the 
more accurate approximation (for clarity the Table includes the number of terms 
needed for the three different expansions). 
With respect to the comparison between the different expansions, we observe
that the computation of the first expansion in terms
of Bessel functions is, as expected, extremely efficient in its range of validity. 
On the other hand, the expansion (\ref{eq:LagBeszer02}) in terms
of Bessel functions is slightly more efficient than the Airy expansion.

\renewcommand{\arraystretch}{1.2}
\begin{table}
\caption{Laguerre expansions, $\alpha$ small: average CPU times per node. 
In this table $\alpha=1/4$. 
The algorithms have been implemented in Fortran 90.
The nodes
are computed with $15-16$ digits accuracy.
\label{tab:tab03}}
$$
\begin{array}{rcc}
{\rm Expansion} & 2\  {\rm terms} \ & \ge 2\ {\rm terms} \quad\\
\hline
  & n=100 & \\
\hline
{\rm Bessel 1}  &   0.03\,\mu s  &   0.16\,\mu s\mbox{ (5 terms)}  \\
{\rm Bessel 2}  &   0.55\,\mu s  &   0.78\,\mu s \mbox{ (4 terms)}   \\
{\rm Airy}  &   0.75\,\mu s  &   1.1\,\mu s  \mbox{ (3 terms)}\\
\hline
& n=10000 & \\
\hline
{\rm Bessel 1}  &   0.03\,\mu s  &   0.07\,\mu s\mbox{ (3 terms)}  \\
{\rm Bessel 2} &   0.53\,\mu s  &   0.53\,\mu s \mbox{ (2 terms)}   \\
{\rm Airy}  &   0.84\,\mu s  &   0.84\,\mu s \mbox{ (2 terms)}  \\
\end{array}
$$
\end{table}
\renewcommand{\arraystretch}{1.0}

 As in the case of the Gauss--Hermite quadrature, overflow/underflow limitations
 in the computation of the weights can be eliminated by balancing the large
 terms as a function of $n$ in the expressions and by scaling out the dependence
 on the weights. A first estimation of the weights as $n\rightarrow \infty$ is
 given by 
 \begin{equation}
 w_i\sim \Frac{\pi}{\sqrt{n}}x_i^{\alpha+1/2}e^{-x_i}.
 \end{equation}
The range of computation of the weights of
Gauss--Laguerre quadrature can be enlarged by simply scaling out the dominant
factor in the asymptotic expansions for Laguerre polynomials. When $\alpha$ is small, this factor is given by $e^{x/2}$. With this, one can define the scaled weights by

\begin{equation}
\label{wghsLag}
\tilde{w}_i =w_i e^{x_i} x_i^{\alpha+1/2}\,.
\end{equation}

These normalized weights do not overflow/underflow as a function of $n$, $\alpha$ and $x_i$. In addition,
 similarly as we did for the Hermite case, we can compute this scaled weights in a numerically stable
way. We notice that the weights (\ref{wgl}) can be written as
\begin{equation}
w_i=\Frac{4\Gamma (n+\alpha+1)}{n! \left[\Frac{d}{dz}L_n^{(\alpha)}(z_i^2)\right]^2},
\end{equation}
where $z=\sqrt{x}$, and therefore $z_i=\sqrt{x_i}$.
Now, in the new variable $z$, the scaled weights can be expressed as 
\begin{equation}
\label{weightsL}
\tilde{w}_i=\Frac{4\Gamma (n+\alpha +1)}{n! (\dot{y}(z_i))^2},
\end{equation}
where the dots mean differentiation with respect to $z$ and
\begin{equation}
y(z)=z^{\alpha+1/2}e^{-z^2/2}L_n^{(\alpha)}(z^2).
\end{equation}
Now, we define $W(z)=4\Gamma (n+\alpha +1)/(n! (\dot{y}(z))^2)$ and with this we have that
$w_i=W(z_i)$, and it is straightforward to check that we have again the desirable property
$\frac{d}{dz}W(z_i)=0$. This means that the computation is well conditioned in the sense that
the error for the weights will be approximately proportional to the square of the error for
the nodes. As a consequence, as we will shown, the weights can can be computed with almost no accuracy loss.

All that is left for computing the nodes is to use the expansions for the Laguerre polynomials in order
to compute $\dot{y}(z)$ by differentiation. In particular, starting from (\ref{eq:lagbess05}) we
have
\begin{equation}
\dot{y}(z)=\left(\frac{\nu}{2}\right)^{\frac{\alpha -1}{2}}\left(\frac{2\nu\zeta \rho^2}{1-\rho}\right)^{1/4}
\left[J_{\alpha} \left(2\nu \sqrt{\zeta}\right)C(\zeta)-\Frac{1}{\sqrt{\zeta}}J_{\alpha+1}
\left(2\nu\sqrt{\zeta})\right) D(\zeta)\right],
\end{equation}
where in this expression $x$ is the variable defined in Section (\ref{sec:LagBesexp}) 
and
\begin{equation}
\begin{array}{l}
C(\zeta)=\left\{\Frac{1}{4(1-\rho)}+\left(\frac12+\alpha\right)\varphi\right\}A+A^{\prime}
-2\nu\varphi B,\\
D(\zeta)=\left\{\Frac{1}{4(1-\rho)}-\left(\frac32+\alpha\right)\varphi\right\}B+B^{\prime}
+2\nu\zeta \varphi A,
\end{array}
\end{equation}
and in these equations prime denotes the derivative with respect to $\rho$.

\begin{figure}
\begin{center}
\hspace*{-1.5cm}
\epsfxsize=15cm \epsfbox{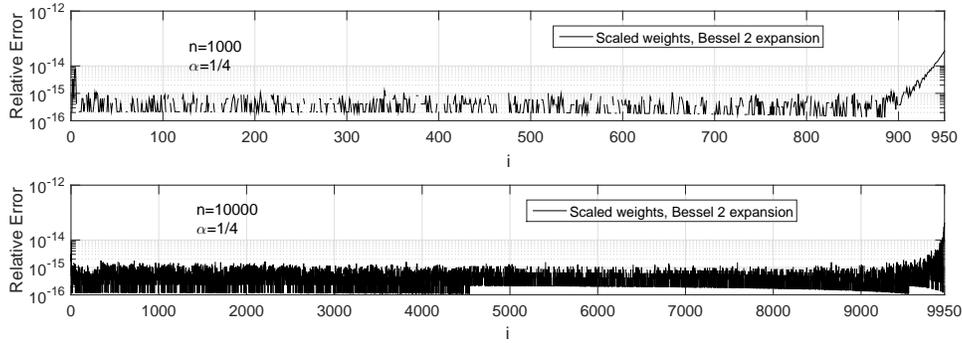}
\caption{ Relative accuracy obtained for the computation of the scaled weights (\ref{scaledw})
for $n=1000,\,10000$ (with $\alpha=1/4$) using the asymptotic expansion for the 
Laguerre polynomials (\ref{eq:lagbess05}) in terms of Bessel functions. 
\label{FigweightsL}}
\end{center}
\end{figure}

Similarly as we did for the Hermite case, we show in Figure \ref{FigweightsL} two examples 
of computation of the scaled weights (\ref{weightsL}) for $n=1000,\,10000$
(with $\alpha=1/4$). 
We use the expansion in terms of Bessel functions 
 (\ref{eq:lagbess05}).  As can be seen, the accuracy for the scaled weights is 
better than $10^{-15}$ in most cases. There is some loss of accuracy for the weights corresponding to the largest nodes
(as discussed, for these values one has to use the expansion for the Laguerre
polynomials in terms of Airy functions).

\subsection{Expansions for large  values of  \protectbold{\alpha}}
\label{lalp}

\subsubsection{An expansion for large  values of  \protectbold{\alpha} and fixed degree \protectbold{n}}\label{sec:Lagherla}

From the well-known limit 
\begin{equation}\label{eq:lagalphaher01}
\lim_{\alpha\to \infty} \alpha^{-n} L_{n}^{(\alpha)}(\alpha t) = \frac{(1-t)^{n}}{n!},
\end{equation}
it follows that the zeros of $L_{n}^{(\alpha)}(\alpha t)$ coalesce at $t=1$ when $\alpha$ is large and $n\ll\alpha$. The limit gives limited information when $t=1$, and in this section we give more details about the behavior of $L_{n}^{(\alpha)}(\alpha t)$ for small values of $\vert t-1\vert$. We consider an asymptotic representation in terms of Hermite polynomials, which has been derived in \cite{Lopez:1999:AOP}. 

We have
\begin{equation}\label{eq:lagalphaher02}
L_{n}^{(\alpha)}(x) =(-1)^n\,z^n\,\sum_{k=0}^n\,\frac{c_k}{z^k}\,
\frac{H_{n-k}(\zeta)}{(n-k)!},  
 \end{equation}
where
\begin{equation}\label{eq:lagalphaher03}
z=\sqrt{{x-(\alpha+1)/2}},\quad \zeta=\frac{x-\alpha-1}{2z}.
\end{equation}
The representation in \eqref{eq:lagalphaher02} holds for $n=0,1,2,\ldots$, and all complex 
values of  $x$ and $\alpha$ and has an asymptotic character for large values
of $|\alpha|+|x|$; the degree $n$ should be fixed. It is not difficult to verify that the limit given in \eqref{eq:lagalphaher01}  follows from \eqref{eq:lagalphaher02}.

The coefficients $c_k$ are defined by
\begin{equation}\label{eq:lagalphaher04}
c_0=1, \quad c_1=c_2=0, \quad c_3= \tfrac13(3x-\alpha-1),\quad 
c_4=\tfrac14(-4x+\alpha+1),
\end{equation}
and  the recursion relation
\begin{equation}\label{eq:lagalphaher05}
kc_k=-2(k-1)c_{k-1}-(k-2)c_{k-2}+(3x-\alpha-1)c_{k-3}+
(2x-\alpha-1)c_{k-4}.
\end{equation}

An approximation of the zeros of $L_{n}^{(\alpha)}(x)$  can be found in \cite{Calogero:1978:ABZ}, and in  \cite{Lopez:1999:AOP} it is shown that  it can be derived from the expansion given in \eqref{eq:lagalphaher02}. Calogero's result is
\begin{equation}\label{eq:lagalphaher06}
\ell_{n,m}= \alpha+\sqrt{{2\alpha}}h_{n,m}+\tfrac13(1+2n+2h_{n,m}^2)+\bigO\left(\alpha^{-\frac12}\right),\quad \alpha\to\infty,
\end{equation}
where $\ell_{n,m}$ and $h_{n,m}$ denote the corresponding zeros of the Laguerre and Hermite polynomials.

For example, with $n=10$ and $\alpha=1000$, the relative error is not larger than $0.85\times 10^{-3}$ (for the first zero). For the fifth and sixth zero the relative errors are about $0.65\times 10^{-4}$.

\subsubsection{An expansion for large  values of  \protectbold{n} and  \protectbold{\alpha}}\label{sec:LagBesla}
In \cite{Temme:1986:LLD} we have given expansions for large $n$ in which $\alpha=\bigO(n)$ is allowed; for a summary see 
\cite{Temme:1990:AEL}. The results  
follow also from uniform expansions of
Whittaker functions obtained by using differential equations; see 
\cite{Dunster:1989:UAE}. These expansions include the $J$-Bessel function, and are valid in the parameter domain where order and argument of the Bessel function are equal, that is, in the turning point domain.  In this section, explicit expressions for the first few coefficients of the expansion
are given.

By using an integral  we can derive the following asymptotic representation
\begin{equation}\label{eq:LagBesla01}
L_n^{(\alpha)}(4\kappa x)=e^{-\kappa A}\chi(b)\left(\frac{b}{2\kappa x}\right)^{\alpha}\frac{\Gamma(n+\alpha+1)}{n!}
\left(J_\alpha(4\kappa b)A(b)-2bJ_\alpha^\prime(4\kappa x)B(b)\right),
\end{equation}
with expansions
\begin{equation}\label{eq:LagBesla02}
A(b)\sim\sum_{k=0}^\infty \frac{A_k(b)}{\kappa^k}, \quad
B(b)\sim\sum_{k=0}^\infty \frac{B_k(b)}{\kappa^k},
\end{equation}
where
\begin{equation}\label{eq:LagBesla03}
\kappa=n+\tfrac12(\alpha+1),\quad \chi(b)=\left(\frac{4b^2-\tau^2}{4x-4x^2-\tau^2}\right)^{\frac14},\quad \tau=\frac{\alpha}{2\kappa}. \end{equation}
We assume that $\tau <1$. The quantity  $b$ is a function of $x$ and follows from the relation
\begin{equation}\label{eq:LagBesla04}
\begin{array}{@{}r@{\;}c@{\;}l@{}}
&&\dsp{2W-2\tau\arctan\frac{W}{\tau}}=\\[8pt]
&&\quad\quad\quad
\dsp{2R-\arcsin\frac{1-2x}{\sqrt{1-\tau^2}}-\tau\arcsin\frac{x-\frac12\tau^2}{x\sqrt{1-\tau^2}}+\tfrac12\pi(1-\tau),}
\end{array}
\end{equation}
where 
\begin{equation}\label{eq:LagBesla05}
R=\tfrac12\sqrt{4x-4x^2-\tau^2}=\sqrt{(x_2-x)(x-x_1)},\quad W=\sqrt{4b^2-\tau^2},
\end{equation}
and
\begin{equation}\label{eq:LagBesla06}
x_1=\tfrac12\left(1-\sqrt{1-\tau^2}\right), \quad 
x_2=\tfrac12\left(1+\sqrt{1-\tau^2}\right).
\end{equation}

The relation in \eqref{eq:LagBesla04} can be used for $x\in[x_1,x_2]$, in which case $b\ge\frac12\tau$. In this interval the zeros of $L_n^{(\alpha)}(4\kappa x)$ occur. For $x$ outside this interval we refer to \cite[\S3.1]{Temme:1990:AEL}.

The first coefficients of the expansions in \eqref{eq:LagBesla02} are 
\begin{equation}\label{eq:LagBesla07}
\begin{array}{@{}r@{\;}c@{\;}l@{}}
A_0(b)&=&1,\quad B_0(b)=0,\\[8pt]
A_1&=&\dsp{\frac{\tau}{24(\tau^2-1)}}, \quad B_1=\dsp{\frac{P R^3+QW^3}{192R^3W^4(\tau^2-1)}},\\[8pt]
&&\quad P= 4(2\tau^2+12b^2)(1-\tau^2),\\[8pt]
&&\quad Q= 2\tau^4-12x^2\tau^2-\tau^2-8x^3+24x^2-6x,\\[8pt]
B_2(b&)=&A_1(b)B_1(b).
\end{array}
\end{equation}

\subsubsection{Expansions of the zeros}\label{sec:LagLargealpha}
A zero  of $L_n^{(\alpha)}(4\kappa x)$ is a zero of $U(b)$ defined by 
\begin{equation}\label{eq:LagBesla08}
U(b)=J_\alpha(4\kappa b)A(b)-2bJ_\alpha^\prime(4\kappa x)B(b),
\end{equation}
where the relation between $b$ and $x$ is given in \eqref{eq:LagBesla04}. We write a zero in terms of $b$ in the form
\begin{equation}\label{eq:LagBesla09}
b=b_0+\eps, \quad b_0=\frac{j_k}{4\kappa}
\end{equation}
where $j_k$ is a zero of the Bessel function $J_\alpha(z)$. We assume for $\eps$ an expansion in the form
\begin{equation}\label{eq:LagBesla10}
\eps\sim\frac{b_1}{\kappa}+\frac{b_2}{\kappa^2}+\frac{b_3}{\kappa^3}+\ldots.
\end{equation}
By expanding $U(b)$ at $b_0$ we have
\begin{equation}\label{eq:LagBesla11}
U(b_0)+\frac{\eps}{1!}U^\prime(b_0)+ \frac{\eps^2}{2!}U^{\prime\prime}(b_0)+\ldots = 0.
\end{equation}
Using the representation of $U(b)$ given in \eqref{eq:LagBesla08},  substituting the expansion of $\eps$, those of $A(b)$ and $B(b)$ given in \eqref{eq:LagBesla02}, and comparing equal powers of $\kappa$, we can obtain the coefficients $b_j$ of \eqref{eq:LagBesla12}.

The first coefficients are
\begin{equation}\label{eq:LagBesla12}
\begin{array}{@{}r@{\;}c@{\;}l@{}}
b_1&=&0,\quad b_2=\frac12 b B_1(b),\quad b_3=\frac12 b\left(B_2(b)-A_1(b)B_1(b)\right)=0, \\[8pt]
b_4&=&\frac1{24}b\left(12B_3(b)-16b^2B_1^3(b)+6bB_1^\prime(b)B_1(b)-12A_2(b)B_1(b)+3B_1^2(b)\right),
\end{array}
\end{equation}
with $b=b_0$ given in \eqref{eq:LagBesla09}.

For example, when we take $n=100$, $\alpha=75$, then we obtain for the first zero $b_0=0.1504907582034649$. We find with  this value for $b$ from \eqref{eq:LagBesla04} a first approximation $x=0.0231157462791716$, with a relative error $2.45\times10^{-5}$. We compute with this $x$ and $b=b_0$ the coefficient $b_2$ and find from $b\sim b_0+b_2/\kappa^2$ the value $b=0.1504905751793771$. Again inverting \eqref{eq:LagBesla04} to find the corresponding $x$-value, we find $x=0.0231156896044437$, now with relative error $3.01618\times10^{-11}$.

\section{Acknowledgements}

The authors thank the referees for their constructive remarks. The authors acknowledge financial support from {\it Ministerio de Econom\'{\i}a y Competitividad}, 
project MTM2015-67142-P (MINECO/FEDER, UE).
NMT thanks CWI, Amsterdam, for scientific support.

\bibliographystyle{plain}    
\bibliography{Gauss} 

\end{document}